\documentclass[aop,MSNbibl,citesort,dvips]{arximspdf}
\usepackage{mathbh}
\usepackage{graphicx}

\doi{10.1214/10-AOP542}
\volume{38}
\issue{6}
\pubyear{2010}
\firstpage{2443}
\lastpage{2485}

\makeatletter

\newtheorem{theorem}{Theorem}
\newtheorem{lem}[theorem]{Lemma}
\newtheorem{cor}[theorem]{Corollary}
\newtheorem{prop}[theorem]{Proposition}

\newproclaim{remark}{Remark}

\renewcommand{\mid}{|}
\newcommand{\ind}{\mathbh{1}}
\newcommand{\eps}{\varepsilon}

\newcommand{\cX}{\mathcal{X}}
\newcommand{\vone}{\mathbf{1}}
\newcommand{\vX}{\mathbf{X}}
\newcommand{\vY}{\mathbf{Y}}
\newcommand{\vZ}{\mathbf{Z}}
\newcommand{\mve}{\mathbf{e}}
\newcommand{\mvv}{\mathbf{v}}
\newcommand{\mvgs}{\bolds{\sigma}}

\newcommand{\bB}{\mathbb{B}}
\newcommand{\bZ}{\mathbb{Z}}

\newcommand{\dR}{\mathbb{R}}
\newcommand{\E}{\mathbb{E}}
\newcommand{\pr}{\mathbb{P}}
\newcommand{\var}{\operatorname{Var}}

\makeatother

\begin{document}
\begin{frontmatter}

\title{Applications of Stein's method for concentration~inequalities}
\runtitle{Stein's method for concentration inequalities}

\begin{aug}
\author[A]{\fnms{Sourav} \snm{Chatterjee}\corref{}\thanksref{t1}\ead[label=e1]{sourav@stat.berkeley.edu}} and
\author[A]{\fnms{Partha S.} \snm{Dey}\ead[label=e2]{partha@stat.berkeley.edu}}
\runauthor{S. Chatterjee and P. S. Dey}
\affiliation{University of California, Berkeley}
\address[A]{Department of Statistics\\
University of California, Berkeley\\
367 Evans Hall \#3860\\
Berkeley, California 94720-3860\\
USA\\
\printead{e1}\\
\phantom{E-mail: }\printead*{e2}} %adresu isvedimo komanda gale!
\end{aug}

\thankstext{t1}{Supported in part by NSF Grant DMS-07-07054
and a Sloan Research Fellowship.}

% HISTORY:
\received{\smonth{6} \syear{2009}}
\revised{\smonth{2} \syear{2010}}

% ABSTRACT
%
\begin{abstract}
Stein's method for concentration inequalities was introduced to prove
concentration of measure in problems involving complex dependencies
such as random permutations and Gibbs measures. In this paper, we
provide some extensions of the theory and three applications: (1) We
obtain a concentration inequality for the magnetization in the
Curie--Weiss model at critical temperature (where it obeys a nonstandard
normalization and super-Gaussian concentration). (2) We derive exact
large deviation asymptotics for the number of triangles in the Erd\H
os--R\' enyi random graph $G(n,p)$ when $p \ge0.31$. Similar results
are derived also for general subgraph counts. (3) We obtain some
interesting concentration inequalities for the Ising model on lattices
that hold at all temperatures.
\end{abstract}

% KEYWORDS
%
\begin{keyword}[class=AMS]
\kwd[Primary ]{60E15}
\kwd{60F10}
\kwd[; secondary ]{60C05}
\kwd{82B44}.
\end{keyword}
\begin{keyword}
\kwd{Stein's method}
\kwd{Gibbs measures}
\kwd{concentration inequality}
\kwd{Ising model}
\kwd{Curie--Weiss model}
\kwd{large deviation}
\kwd{Erd\H{o}s--R\'enyi random graph}
\kwd{exponential random graph}.
\end{keyword}

\pdfkeywords{60E15, 60F10, 60C05, 82B44, Stein's method,
Gibbs measures, concentration inequality, Ising model,
Curie--Weiss model, large deviation, Erdos--Renyi random graph,
exponential random graph}

\end{frontmatter}

%s1 ###
\section{Introduction}\label{sec:intro}
In his seminal 1972 paper \cite{stein72}, Charles Stein introduced a
method for proving central limit theorems with convergence rates for
sums of dependent random variables. This has now come to be known as
\textit{Stein's method}. The technique is primarily used for proving
distributional limit theorems (both Gaussian and non-Gaussian). Stein's
attempts \cite{stein86} at devising a version of the method for large
deviations did not prove fruitful. Some progress for sums of dependent
random variables was made by Rai\v{c} \cite{raic04}. The problem was
finally solved in full generality in \cite{chatterjee05}. A selection
of results and examples from \cite{chatterjee05} appeared in the later
papers \cite{chatterjee07,chatterjee07a}. In this paper, we extend the
theory and work out three further examples. The paper is fully self-contained.

The sections are organized as follows. In Section \ref{sec:results}, we
state the main results, the examples, and some proof sketches. The
complete proofs are in Section \ref{sec:proof}.

%s2 ###
\section{Results and examples}\label{sec:results}
%%%%%%%%%%%%%%%%%%%%%%%%%%%
The following abstract theorem is quoted\break from \cite{chatterjee07}. It
summarizes a collection of results from \cite{chatterjee05}. This is a
generalization of Stein's method of exchangeable pairs to the realm of
concentration inequalities and large deviations.
\begin{theorem}[(\cite{chatterjee07}, Theorem 1.5)]\label{thm:conc}
Let $\cX$ be a separable metric space and suppose $(X,X')$ is an
exchangeable pair of $\cX$-valued random variables. Suppose $f\dvtx\cX\to
\dR$ and $F\dvtx\cX\times\cX\to\dR$ are square-integrable functions such
that $F$ is antisymmetric [i.e., $F(X,X')=-F(X',X)$ a.s.], and
$\E(F(X,X')\mid X) = f(X)$ a.s. Let
\[
\Delta(X) := \tfrac{1}{2}\E\bigl(\bigl|\bigl(f(X)-f(X')\bigr)F(X,X')\bigr| |
X\bigr).
\]
Then $\E(f(X)) = 0$, and the following concentration results hold for $f(X)$:
\begin{longlist}
\item If $\E(\Delta(X)) <\infty$, then $\var(f(X)) = \frac
{1}{2}\E((f(X)-f(X'))F(X,X'))$.
\item Assume that $\E(e^{\theta f(X)} |F(X,X')|) < \infty$ for
all $\theta$. If there exists nonnegative constants $B$ and $C$ such
that $\Delta(X) \le B f(X) + C$ almost surely, then for any $t\ge0$,
\[
\pr\{f(X) \ge t\}\le\exp\biggl(-\frac{t^2}{2C + 2Bt} \biggr)
\quad\mbox
{and}\quad \pr\{f(X) \le-t\}\le\exp\biggl(-\frac{t^2}{2C} \biggr).
\]
\item For any positive integer $k$, we have the following
exchangeable pairs version of the Burkholder--Davis--Gundy inequality:
\[
\E(f(X)^{2k}) \le(2k-1)^k \E(\Delta(X)^k).
\]
\end{longlist}
\end{theorem}

Note that the finiteness of the exponential moment for all $\theta$
ensures that the tail bounds hold for all $t$. If it is finite only in
a neighborhood of zero, the tail bounds will hold for $t$ less than a threshold.

One of the contributions of the present paper is the following
generalization of the above result for non-Gaussian tail behavior. We
apply it to obtain a concentration inequality with the correct tail
behavior in the Curie--Weiss model at criticality.
\begin{theorem}\label{thm:nong1}
Suppose $(X,X')$ is an exchangeable pair of random variables. Let
$F(X,X'), f(X)$ and $\Delta(X)$ be as in Theorem \ref{thm:conc}. Suppose
that we have
\[
\Delta(X) \le\psi(f(X)) \qquad\mbox{almost surely}
\]
for some nonnegative symmetric function $\psi$ on $\dR$. Assume that
$\psi$ is nondecreasing and twice continuously differentiable in
$(0,\infty)$ with
%
%e2 ###
%e1 ###
\begin{equation}
\label{eq:al}
\alpha := \sup_{x>0}x\psi'(x)/\psi(x)<2
\end{equation}
and
\begin{equation}
\label{eq:gd}
\delta := \sup_{x>0}x\psi''(x)/\psi(x)<\infty.
\end{equation}
Assume that $\E(|f(X)|^{k}) < \infty$ for all positive integer $k\ge
1$. Then
for any $t\ge0$, we have
\[
\pr\bigl(\vert f(X) \vert> t\bigr) \le c\exp\biggl(-\frac{t^2}{2\psi(t)} \biggr)
\]
for some constant $c$ depending only on $\alpha,\delta$. Moreover, if
$\psi$ is only once differentiable with $\alpha<2$ as in (\ref{eq:al}),
then the tail inequality holds with exponent $t^2/4\psi(t)$.
\end{theorem}

An immediate corollary of Theorem \ref{thm:nong1} is the following.
\begin{cor}\label{cor:nong}
Suppose $(X,X')$ is an exchangeable pair of random variables. Let
$F(X,X'), f(X)$ and $\Delta(X)$ be as in Theorem \ref{thm:conc}. Suppose
that for some real number $\alpha\in(0,2)$ we have
\[
\Delta(X) \le B \vert f(X) \vert^\alpha+ C \qquad\mbox{almost surely},
\]
where $B>0,C\ge0$ are constants. Assume that $\E(|f(X)|^{k}) < \infty$
for all positive integer $k\ge1$. Then for any $t\ge0$ we have
\[
\pr\bigl(\vert f(X) \vert> t\bigr) \le c_{\alpha}\exp\biggl(-\frac{1}{2}\cdot
\frac
{t^2}{Bt^\alpha+C} \biggr)
\]
for some constant $c_{\alpha}$ depending only on $\alpha$.
\end{cor}

The result in Theorem \ref{thm:nong1} states that the tail behavior of
$f(X)$ is essentially given by the behavior of $f(X)^{2}/\Delta(X)$.
Condition (\ref{eq:al}) implies that $\psi(x)< \psi(1)(1+x^2)$ for all
$x\in\dR$. Moreover, the constant $c_{\alpha}$ appearing in
Theorem \ref
{thm:nong1} can be written down explicitly but we did not attempt to
optimize the constant.
The proof of Theorem \ref{thm:nong1} is along the same lines as Theorem
\ref{thm:conc}, but somewhat more involved. Deferring the proof to
Section \ref{sec:proof}, let us move on to examples.
%%%%%%%%%%%%%%%%%%%%%%%%%%%

%s2.1 ###
\subsection{Example: Curie--Weiss model at criticality}\label{subsec:curie}
%%%%%%%%%%%%%%%%%%%%%%%%%%%

The ``Curie--Weiss model of ferromagnetic interaction'' at inverse
temperature $\beta$ and zero external field is given by the following
Gibbs measure on $\{+1,-1\}^{n}$. For a typical configuration $\mvgs
=(\sigma_{1},\sigma_{2},\ldots,\sigma_{n})\in\{+1,-1\}^{n}$, the
probability of
$\mvgs$ is given by
\[
\mu_{\beta}(\{\mvgs\}) := Z_{\beta}^{-1} \exp\biggl( \frac{\beta
}{n}\sum_{i<j}
\sigma_{i}\sigma_{j} \biggr),
\]
where $Z_{\beta}=Z_{\beta}(n)$ is the normalizing constant.
It is well known that the Curie--Weiss model shows a phase transition at
$\beta_{c}=1$. For $\beta<\beta_{c}$, the magnetization $m(\mvgs
):=\frac
{1}{n}\sum_{i=1}^{n}\sigma_{i}$ is concentrated at $0$ but for $\beta
>\beta
_{c}$ the magnetization is concentrated on the set $\{-x^{*},x^{*}\}$
where $x^{*}>0$ is the largest solution of the equation $x=\tanh(\beta
x)$. In fact, using concentration inequalities for exchangeable pairs
it was proved in \cite{chatterjee05} (Proposition 1.3) that for all
$\beta
\ge0,h\in\dR, n\ge1, t\ge0$ we have
\[
\pr\biggl(|m-\tanh(\beta m+h)|\ge\frac{\beta}{n} + \frac{t}{\sqrt
{n}}
\biggr) \le2\exp\biggl(-\frac{t^{2}}{4(1+\beta)} \biggr),
\]
where $h$ is the external field, which is zero in our case. Although a
lot is known about this model (see Ellis \cite{ellis85}, Section IV.4,
for a survey), the above result---to the best of our knowledge---is
the first rigorously proven concentration inequality that holds at all
temperatures. (See also \cite{chazottes06} for some related results.)

Incidentally, the above result shows that when $\beta< 1$, the
magnetization is at most of order $n^{-1/2}$. It is known that at the
critical temperature the magnetization $m(\mvgs)$ shows a non-Gaussian
behavior and is of order $n^{-1/4}$. In fact, at $\beta=1$ as $n\to
\infty
$, $n^{1/4}m(\mvgs)$ converges to the probability distribution on $\dR$
having density proportional to $\exp(-t^{4}/12)$. This limit theorem
was first proved by Simon and Griffiths \cite{simon73} and error bounds
were obtained recently \cite{chatterjeeshao09,eichelsbacherlowe09}.
The following concentration inequality, derived using Theorem \ref
{thm:nong1}, fills the gap in the tail bound at the critical point.
\begin{prop}\label{pr:crit}
Suppose $\mvgs$ is drawn from the Curie--Weiss model at the critical
temperature $\beta=1$. Then, for any $n\ge1$ and $ t\ge0$, the
magnetization satisfies
\[
\pr\bigl(n^{1/4}|m(\mvgs)| \ge t\bigr) \le2e^{-ct^4},
\]
where $c>0$ is an absolute constant.
\end{prop}

Here, we may remark that such a concentration inequality probably
cannot be obtained by application of standard off-the-shelf results
(e.g., those surveyed in Ledoux \cite{ledoux01}, the famous results of
Talagrand \cite{talagrand95} or the recent breakthroughs of Boucheron,
Lugosi and Massart \cite{blm03}), because they generally give Gaussian
or exponential tail bounds. There are several recent remarkable results
giving tail bounds different from exponential and Gaussian. The papers
\cite{bcr06,lo00,ggm05,chazottes06} deal with tails between exponential and
Gaussian and \cite{bcr05,bob07} deal with subexponential tails. Also in
\cite{bl00,gozlan07,gozlan09}, the authors deal with tails (possibly)
larger than Gaussian. However, it seems that none of the techniques
given in these references would lead to the result of Proposition \ref
{pr:crit}.

It is possible to derive a similar tail bound using the asymptotic
results of Martin-L\"{o}f \cite{martinlof82} about the partition
function $Z_{\beta}(n)$ (see also Bolthausen \cite{bolthausen87}). An
application of their results gives that
\[
\sum_{\mvgs\in\{-1,+1\}^{n}}e^{{n}/{2}m(\mvgs)^{2}+n\theta
m(\mvgs
)^{4}}\simeq\frac{2^{n+1}\Gamma(5/4)}{\sqrt{2\pi}} \biggl(\frac
{12n}{1-12\theta} \biggr)^{1/4}
\]
for $\theta<1/12$ in the sense that the ratio of the two sides
converges to one as $n$ goes to infinity and from here the tail bound
follows easily (without an explicit constant).
However, this approach depends on a precise estimate of the partition
function [e.g., large deviation estimates or finding the
limiting free energy $\lim n^{-1}\log Z_{\beta}(n)$ are not enough] and
this precise estimate is hard to prove. Our method, on the other hand,
depends only on simple properties of the Gibbs measure and is not tied
specifically to the Curie--Weiss model.

The idea used in the proof of Proposition \ref{pr:crit} can be used to
prove a tail inequality that holds for all $0\le\beta\le1$. We state the
result below without proof. Note that the inequality gives the correct
tail bound for all $0\le\beta\le1$.
\begin{prop}
Suppose $\mvgs$ is drawn from the Curie--Weiss model at inverse
temperature $\beta$ where $0\le\beta\le1$. Then, for any $n\ge1$
and $
t\ge0$ the magnetization satisfies
\[
\pr\bigl(3(1-\beta)m(\mvgs)^{2}+\beta^{3}m(\mvgs)^{4} \ge t\bigr) \le2e^{-nt/160}.
\]
\end{prop}

It is possible to derive similar non-Gaussian tail inequalities for
general Curie--Weiss models at the critical temperature. We briefly
discuss the general case below. Let $\rho$ be a symmetric probability
measure on $\dR$ with $\int x^{2} \,d\rho(x)=1$ and $\int\exp(\beta
x^{2}/2) \,d\rho(x)<\infty$ for all $\beta\ge0$.
The general Curie--Weiss model CW$(\rho)$ at inverse temperature $\beta$
is defined as the array of spin random variables $\vX
=(X_{1},X_{2},\ldots,X_{n})$ with joint distribution
%
%e3 ###
\begin{equation}
d\nu_{n}(\mathbf{x})=Z_{n}^{-1}\exp\biggl(\frac{\beta}{2n}
(x_{1}+x_{2}+\cdots
+x_{n} )^{2} \biggr)\prod_{i=1}^{n}d\rho(x_{i})
\end{equation}
for $\mathbf{x}=(x_{1},x_{2},\ldots,x_{n})\in\dR^{n}$
where
\[
Z_{n}=\int\exp\biggl(\frac{\beta}{2n}(x_{1}+x_{2}+\cdots
+x_{n})^{2}
\biggr)\prod_{i=1}^{n}d\rho(x_{i})
\]
is the normalizing constant. The magnetization $m(\mathbf{x})$ is
defined as
usual by $m(\mathbf{x})=n^{-1}\sum_{i=1}^{n}x_{i}$. Here, we will
consider the
case when $\rho$ satisfies the following two conditions:
\begin{enumerate}[(B)]
\item[(A)] $\rho$ has compact support, that is, $\rho([-L,L])=1$ for some
$L<\infty$.
\item[(B)] The equation $h'(s)=0$ has a unique root at $s=0$ where
\[
h(s):=\frac{s^{2}}{2}-\log\int\exp(sx) \,d\rho(x)\qquad \mbox{for }
s\in\dR.
\]
\end{enumerate}
The second condition says that $h(\cdot)$ has a unique global minima at
$s=0$ and $|h'(s)|>0$ for $|s|>0$.
The behavior of this model is quite similar to the classical
Curie--Weiss model and there is a phase transition at $\beta=1$. For
$\beta
<1$, $m(\vX)$ is concentrated around zero while for $\beta>1, m(\vX
)$ is
bounded away from zero a.s. (see Ellis and Newman \cite{ellis78,ellis78a}).
We will prove the following concentration result.
\begin{prop}\label{pr:critgen}
Suppose $\vX\sim\nu_{n}$ at the critical temperature $\beta=1$ where
$\rho
$ satisfies conditions \textup{(A)} and \textup{(B)}. Let $k$ be such that
$h^{(i)}(0)=0$ for $0\le i<2k$ and $h^{(2k)}(0)\neq0$, where
\[
h(s):=\frac{s^{2}}{2}-\log\int\exp(sx) \,d\rho(x) \qquad\mbox{for }
s\in\dR,
\]
and $h^{(i)}$ is the $i$th derivative of $h$. Then, $k>1$ and for any
$n\ge1$ and $ t\ge0$ the magnetization satisfies
\[
\pr\bigl(n^{1/2k}|m(\vX)| \ge t\bigr) \le2e^{-ct^{2k}},
\]
where $c>0$ is an absolute constant depending only on $\rho$.
\end{prop}

Here, we mention that in Ellis and Newman \cite{ellis78}, convergence
results were proved for the magnetization in CW$(\rho)$ model under
optimal condition on $\rho$. Under our assumption, their result says
that $n^{1/2k}m(\vX)$ converges weakly to a distribution having density
proportional to $\exp(-\lambda x^{2k}/(2k)!)$ where $\lambda:=h^{(2k)}(0)$.
Hence, the tail bound gives the correct convergence rate.

Let us now give a brief sketch of the proof of Proposition \ref{pr:crit}.
Suppose $\mvgs$ is drawn from the Curie--Weiss model at the critical
temperature. We construct $\mvgs'$ by taking one step in the heat-bath
Glauber dynamics: a coordinate $I$ is chosen uniformly at random, and
$\sigma_{I}$ is replace by $\sigma'_{I}$ drawn from the conditional
distribution of the $I$th coordinate given $\{\sigma_{j}\dvtx j\neq I\}$. Let
\[
F(\mvgs,\mvgs') := \sum_{i=1}^{n}(\sigma_{i}-\sigma'_{i})=\sigma
_{I}-\sigma'_{I}.
\]
For each $i=1,2,\ldots,n$, define
$
m_{i}=m_{i}(\mvgs) = n^{-1}\sum_{j\neq i} \sigma_{j}.
$
An easy computaion gives that $\E(\sigma_{i}| \{\sigma_{j}, j\neq i\}
) =
\tanh( m_{i})$ for all $i$ and so we have
\[
f(\mvgs):= \E(F(\mvgs,\mvgs')|\mvgs) = m - \frac{1}{n} \sum_{i=1}^{n}
\tanh(m_{i})=\frac{m}{n} + \frac{1}{n} \sum_{i=1}^{n} g(m_{i}),
\]
where $g(x):=x-\tanh(x)$. Note that $|m_i - m|\le1/n$, and hence
$f(\mvgs) = m - \tanh m + O(1/n)$. A simple analytical argument using
the fact that, for $x\approx0$, $x- \tanh x = x^3/3 + O(x^5)$ then gives
\[
\Delta(\mvgs) \le\frac{6}{n}|f(\mvgs)|^{2/3}+\frac{12}{n^{5/3}}
\]
and using Corollary \ref{cor:nong} with $\alpha=2/3, B= 6/n$ and
$C=12/n^{5/3}$ we have
\[
\pr(|m-\tanh m| \ge t+n^{-1}) \le\pr\bigl(|f(\mvgs)| \ge t\bigr) \le2e^{- cnt^{4/3}}
\]
for all $t\ge0$ for some constant $c>0$. It is easy to see that this
implies the result. The critical observation, of course, is that $x -
\tanh(\beta x) = O(x^3)$ for $\beta= 1$, which is not true for $\beta
\ne1$.

%%%%%%%%%%%%%%%%%%%%%%%%%%%%%%%%%%%%%%%%%%%%%%%%
%s2.2 ###
\subsection{Example: Triangles in Erd\H{o}s--R\'enyi graphs}

\begin{figure}

\includegraphics{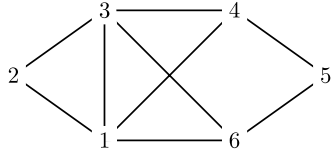}

\caption{A graph on $6$ vertices with $8$ edges and $3$ triangles. The
triangles being $(1,2,3), (1,3,4)$ and $(1,3,6)$.}
\label{fig:Tn}
\end{figure}

Consider the Erd\H{o}s--R\'enyi random graph model $G(n,p)$ which is
defined as follows. The vertex set is $[n]:=\{1,2,\ldots,n\}$ and each
edge $(i,j), 1\le i< j\le n$, is present with probability $p$ and not
present with probability $1-p$ independently of each other. For any
three distinct vertex $i<j<k$ in $[n]$ we say that the triple $(i,j,k)$
forms a triangle in the graph $G(n,p)$ if all the three edges
$(i,j),(j,k),(i,k)$ are present in $G(n,p)$ (see Figure \ref{fig:Tn}).
Let $T_{n}$ be the number of triangles in $G(n,p)$, that is,
\[
T_{n}:=\sum_{1\le i<j<k\le n}\vone\{(i,j,k)\mbox{ forms a triangle
in }
G(n,p)\}.
\]

%f1 ###

Let us define the function $I(\cdot,\cdot)$ on $(0,1)\times(0,1)$ as
%
%e4 ###
\begin{equation}\label{kl}
I(r,s) := r\log\frac{r}{s} + (1-r)\log\frac{1-r}{1-s}.
\end{equation}
Note that $I(r,s)$ is the Kullback--Leibler divergence of the measure
$\nu_{s}$ from $\nu_{r}$ and also the relative entropy of $\nu_{r}$
w.r.t. $\nu_{s}$ where $\nu_{p}$ is the Bernoulli$(p)$ measure. We have
the following result about the large deviation rate function for the
number of triangles in $G(n,p)$.
\begin{theorem}\label{largedev}
Let $T_n$ be the number of triangles in $G(n,p)$, where $p > p_{0}$
where $p_{0}=2/(2+e^{3/2})\approx0.31$. Then for any $r\in(p,1]$,
%
%e5 ###
\begin{equation}\label{eq:trtail}
\pr\left(T_n \ge\pmatrix{n \cr3} r^3 \right) = \exp\biggl(-\frac{n^2
I(r,p)}{2}\bigl(1+O(n^{-1/2})\bigr) \biggr).
\end{equation}
Moreover, even if $p \le p_{0}$, there exist $p',p''$ such that
$p<p'\le p''<1$ and the same result holds for all $r\in(p,p')\cup(p'',1]$.
For all $p$ and $r$ in the above domains, we also have the more precise estimate
%
%e6 ###
\begin{equation}\label{eq:trbound}\quad
\pr\left( \left|T_n- \pmatrix{n \cr3} r^3 \right| \le
C(p,r)n^{5/2} \right)
= \exp\biggl(-\frac{n^2 I(r,p)}{2}\bigl(1+O(n^{-1/2})\bigr) \biggr),
\end{equation}
where $C(p,r)$ is a constant depending on $p$ and $r$.
\end{theorem}

\begin{figure}

\includegraphics{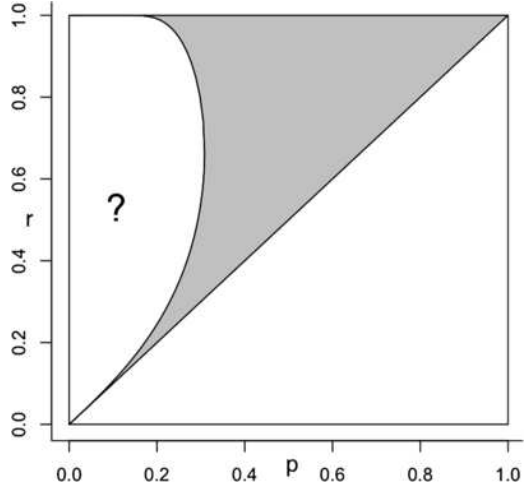}

\caption{The set (colored in gray) of $(p,r), r\ge p$, for which we are
able to show that the large deviation result holds.}
\label{fig:pr}
\end{figure}

The behavior of the upper tail of subgraph counts in $G(n,p)$ is a
problem of great interest in the theory of random graphs (see
\cite{bollobas85,jansonetal00,jansonrucinski02,vu01,kimvu04}, and references
contained therein). The best upper bounds to date were obtained by Kim
and Vu \cite{kimvu04} (triangles) and Janson, Oleszkiewicz, and Ruci\'
nski \cite{jansonolesz04} (general subgraph counts). For triangles, the
results of these papers essentially state that for a fixed $\epsilon>0$,
\[
\exp\bigl(-\Theta\bigl(n^2 p^2\log(1/p)\bigr)\bigr) \le\pr\bigl(T_n \ge\E(T_n) + \epsilon n^3
p^3\bigr) \le\exp(-\Theta(n^2 p^2)).
\]
%
%Lower tails are precisely known (see \cite{jansonetal00}, Chap. 2).
Clearly, our result gives a lot more in the situations where it works
(see Figure \ref{fig:pr}). The method of proof can be easily extended
to prove similar results for general subgraph counts and are discussed
in Section \ref{subsec:small}. However, there is an obvious
incompleteness in Theorem \ref{largedev} (and also for general
subgraphs counts), namely, that it does not work for all $(p,r)$.

%
%f2 ###

In this context, we should mention that another paper on large
deviations for subgraph counts by Bolthausen, Comets and Dembo
\cite{bolthausenetal09} is in preparation. As of now, to the best of our
knowledge, the authors of \cite{bolthausenetal09} have only looked at
subgraphs that do not complete loops, like $2$-stars. Another related
article is the one by D\"oring and Eichelsbacher
\cite{doringeichelsbacher09}, who obtain moderate deviations for a class of
graph-related objects, including triangles.

Unlike the previous two examples, Theorem \ref{largedev} is far from
being a direct consequence of any of our abstract results. Therefore,
let us give a sketch of the proof, which involves a new idea.

The first step is standard: consider tilted measures. However, the
appropriate tilted measure in this case leads to what is known as an
``exponential random graph,'' a little studied object in the rigorous
literature. Exponential random graphs have become popular in the
statistical physics and network communities in recent years (see the
survey of Park and Newman \cite{parknewman04}). The only rigorous work
we are aware of is the recent paper of Bhamidi et al.
\cite{bhamidi08}, who look at convergence rates of Markov chains that
generate such graphs.

We will not go into the general definition or properties of exponential
random graphs. Let us only define the model we need for our purpose.

Fix two numbers $\beta\ge0$ and $h\in\dR$. Let $\Omega=\{0,1\}
^{n\choose2}$ be the space of all tuples like $\mathbf
{x}=(x_{ij})_{1\le
i<j\le n}$, where $x_{ij}\in\{0,1\}$ for each $i,j$. Let $\vX=
(X_{ij})_{1\le i<j\le n}$ be a random element of $\Omega$ following the
probability measure proportional to $e^{H(\mathbf{x})}$, where $H$ is
the Hamiltonian:
\[
H(\mathbf{x}) = \frac{\beta}{n}\sum_{1\le i<j <k\le n}
x_{ij}x_{jk} x_{ik}
+ h\sum_{1\le i<j\le n} x_{ij}.
\]
Note that any element of $\Omega$ naturally defines an undirected graph
on a set of $n$ vertices. For each $\mathbf{x}\in\Omega$, let
$T(\mathbf{x})=\sum
_{i<j<k}x_{ij}x_{jk}x_{ik}$ denote the number of triangles in the graph
defined by $\mathbf{x}$, and let $E(\mathbf{x})=\sum_{i<j}x_{ij}$
denote the number
of edges. Then the above Hamiltonian is nothing but
\[
\frac{\beta T(\mathbf{x})}{n}+ hE(\mathbf{x}).
\]
For notational convenience, we will assume that $x_{ij}=x_{ji}$.
Let $Z_n(\beta,h)$ be the corresponding partition function, that is,
\[
Z_n(\beta,h) = \sum_{\mathbf{x}\in\Omega} e^{H(\mathbf{x})}.
\]
Note that $\beta=0$ corresponds to the Erd\H{o}s--R\'{e}nyi random graph
with $p=e^{h}/(1+e^{h})$.
The following theorem ``solves'' this model in a ``high temperature
region.'' Once this solution is known, the computation of the large
deviation rate function is just one step away.
\begin{theorem}[(Free energy in high temperature regime)]\label{pressure}
Suppose we have $\beta\ge0$, $h\in\dR$, and $Z_n(\beta,h)$ defined
as above. Define a function $\varphi\dvtx[0,1]\to\dR$ as
\[
\varphi(x) = \frac{e^{\beta x + h}}{1+e^{\beta x + h}}.
\]
Suppose $\beta$ and $h$ are such that the equation $u = \varphi(u)^2$
has a unique solution $u^*$ in $[0,1]$ and $2\varphi(u^{*})\varphi
'(u^{*})<1$. Then
\[
\lim_{n\to\infty}\frac{\log Z_n(\beta,h)}{n^2} = -\frac
{1}{2}I(\varphi
(u^*),\varphi(0)) - \frac{1}{2} \log\bigl(1-\varphi(0)\bigr) +\frac{\beta
\varphi
(u^*)^3}{6},
\]
where $I(\cdot,\cdot)$ is the function defined in (\ref{kl}). Moreover,
there exists a constant $K(\beta,h)$ that depends only on $\beta$ and
$h$ (and not on $n$) such that difference between $n^{-2}\log Z_n(\beta
,h)$ and the limit is bounded by $K(\beta,h)n^{-1/2}$ for all $n$.
\end{theorem}

Incidentally, the above solution was obtained using physical heuristics
by Park and Newman \cite{parknewman05} in 2005. Here, we mention that,
in fact, the following result is always true.
\begin{lem}\label{lem:lbd}
For any $\beta\ge0, h\in\dR$, we have
%
%e7 ###
\begin{eqnarray}\label{eq:lbd}
&&\liminf_{n\to\infty}\frac{\log Z_n(\beta,h)}{n^2} \nonumber\\
&&\qquad\ge \sup_{r\in
(0,1)} \biggl\{-\frac{1}{2}I(r,\varphi(0)) - \frac{1}{2} \log
\bigl(1-\varphi
(0)\bigr) +\frac{\beta r^3}{6} \biggr\}\\
&&\qquad= \sup_{u\dvtx\varphi(u)^{2}=u} \biggl\{-\frac{1}{2}I(\varphi
(u),\varphi
(0)) - \frac{1}{2} \log\bigl(1-\varphi(0)\bigr) +\frac{\beta\varphi
(u)^3}{6} \biggr\}.\nonumber
\end{eqnarray}
\end{lem}

We will characterize the set of $\beta,h$ for which the conditions in
Theorem \ref{pressure} hold in Lemma \ref{lem:pcond}. First of all,
note that the appearance of the function $\varphi(u)^{2}-u$ is not
magical. For each $i < j$, define
\[
L_{ij} = \frac{1}{n}\sum_{k\notin\{i,j\}} X_{ik} X_{jk}.
\]
This is the number of ``wedges'' or $2$-stars in the graph that have the
edge $ij$ as base. The key idea is to use Theorem \ref{thm:conc} to
show that these quantities approximately satisfy the following set of
``mean field equations'':
%
%e8 ###
\begin{equation}\label{meanfieldeqs}
L_{ij} \simeq\frac{1}{n}\sum_{k\notin\{i,j\}} \varphi
(L_{ik})\varphi
(L_{jk}) \qquad\mbox{for all } i<j.
\end{equation}
(The idea of using Theorem \ref{thm:conc} to prove mean field equations
was initially developed in Section 3.4 of \cite{chatterjee05}.) The
following lemma makes this notion precise. Later, we will show that
under the conditions of Theorem \ref{pressure}, this system has a
unique solution.
\begin{lem}[(Mean field equations)]\label{lem:mf}
Let $\varphi$ be defined as in Theorem \ref{pressure}. Then for any
$1\le i<j\le n$, we have
\[
\pr\biggl(
\sqrt{n} \biggl|L_{ij} - \frac{1}{n}\sum_{k\notin\{i,j\}}
\varphi
(L_{ik})\varphi(L_{jk}) \biggr|
\ge t \biggr) \le2\exp\biggl(- \frac{t^{2}}{8(1+\beta)} \biggr)
\]
for all $t\ge8\beta/n$.
In particular, we have
%
%e9 ###
\begin{equation}\label{eq:lcon}
\E\biggl|L_{ij} - \frac{1}{n}\sum_{k\notin\{i,j\}} \varphi
(L_{ik})\varphi(L_{jk}) \biggr| \le\frac{C(1+\beta)^{1/2}}{n^{1/2}},
\end{equation}
where $C$ is a universal constant.
\end{lem}

In fact, one would expect that $L_{ij}\simeq u^{*}$ for all $i<j$, if
the equation
%
%e10 ###
\begin{equation}\label{eq:sol}
\psi(u):=\varphi(u)^{2}-u=0
\end{equation}
has a unique solution $u^{*}$ in $[0,1]$. The intuition behind is as
follows. Define $L_{\max} = \max_{i,j} L_{ij}$ and $L_{\min}=\min
_{i,j} L_{ij}$. It is easy to see that $\varphi$ is an increasing
function. Hence, from the mean-field equations (\ref{meanfieldeqs}), we
have $L_{\max}\le\varphi(L_{\max})^{2} + o(1)$ or $\psi(L_{\max})
\ge o(1)$. But $\psi(u)\ge0$ iff $u\le u^{*}$. Hence, $L_{\max}\le
u^{*}+o(1)$. Similarly, we have $L_{\min}\ge u^{*}-o(1)$ and thus all
$L_{ij}\simeq u^{*}$. Lemma \ref{mean2} formalizes this idea. Here, we
mention that one can easily check that equation (\ref{eq:sol}) has at
most three solutions. Moreover, $\psi(0)>0>\psi(1)$ implies that
$\psi
'(u^{*})\le0$ or $2\varphi(u^{*})\varphi'(u^{*})\le1$ if $u^{*}$ is
the unique solution to (\ref{eq:sol}).
\begin{lem}\label{mean2}
Let $u^*$ be the unique solution of the equation $u=\varphi(u)^{2}$.
Assume that $2\varphi(u^{*})\varphi'(u^{*})<1$. Then for each $1\le
i<j\le n$, we have
\[
\E|L_{ij} - u^*| \le\frac{K(\beta,h)}{n^{1/2}},
\]
where $K(\beta,h)$ is a constant depending only on $\beta,h$.
Moreover, if $2\varphi(u^{*})\varphi'(u^{*}) =1$ then we have
\[
\E|L_{ij} - u^*| \le\frac{K(\beta,h)}{n^{1/6}} \qquad\mbox{for all
}1\le
i<j\le n.
\]
\end{lem}

Now observe that the Hamiltonian $H(\vX)$ can be written as
\[
H(\vX) = \frac{\beta}{6} \sum_{1\le i<j\le n} X_{ij} L_{ij} + h\sum
_{1\le i<j\le n} X_{ij}.
\]
The idea then is the following: once we know that the conclusion of
Lemma \ref{mean2} holds, each $L_{ij}$ in the above Hamiltonian can be
replaced by $u^*$, which results in a model where the coordinates are
independent. The resulting probability measure is presumably quite
different from the original measure, but somehow the partition
functions remain comparable.

\begin{figure}[b]

\includegraphics{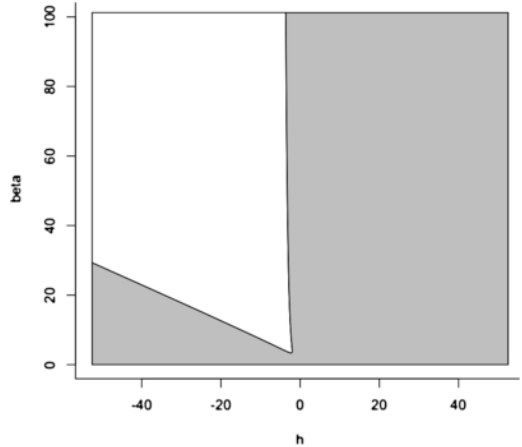}

\caption{The set $S$ (colored in gray) of $(h,\beta)$ for which the
conditions of Theorem \protect\ref{pressure} hold.}
\label{fig:unisol}
\end{figure}

The following lemma (Lemma \ref{lem:pcond}) characterizes the region
$S\in\dR\times[0,\infty)$ such that the equation $u = \varphi(u)^2$
has a unique solution $u^{*}$ in $[0,1]$ and $2\varphi(u^{*})\varphi
'(u^{*})< 1$ for $(h,\beta)\in S$ (see Figure \ref{fig:unisol}).

%
%f3 ###

Let $h_{0}=\log2-\frac{3}{2}<0$. For $h<h_{0}$ there exist exactly two
solutions $0<a_{*}=a_{*}(h)<1/2<a^{*}=a^{*}(h)<\infty$ to the equation
\[
\log x + \frac{1+x}{2x}+h=0.
\]
Define $a_{*}(h)=a^{*}(h)=1/2$ for $h=h_{0}$ and
%
%e11 ###
\begin{equation}\label{eq:bstar}
\beta_{*}(h)= \frac{(1+a_{*})^{3}}{2a_{*}} \quad\mbox{and}\quad \beta^{*}(h)=
\frac
{(1+a^{*})^{3}}{2a^{*}}
\end{equation}
for $h\le h_{0}$.
\begin{lem}[(Characterization of high temperature regime)] \label{lem:pcond}
Let $S$ be the set of pairs $(h,\beta)$ for which the function $\psi
(u):=\varphi(u)^{2}-u$ has a unique root $u^{*}$ in $[0,1]$ and
$2\varphi(u^{*})\varphi'(u^{*})< 1$ where
$
\varphi(u):= e^{\beta u +h}/(1+e^{\beta u +h}).
$
Then we have
\[
S^{c}=\{(h,\beta)\dvtx h\le h_{0} \mbox{ and } \beta_{*}(h) \le\beta
\le\beta
^{*}(h)\},
\]
where $\beta^{*},\beta_{*}$ are as given in equation (\ref{eq:bstar}).
In particular, $(h,\beta)\in S$ if $\beta\le(3/2)^{3}$ or $h>h_{0}$.
\end{lem}
%
% for $$t_{1}=\frac{\sqrt{ t^2(t+3)^2 +4t} -t(t+3)}{2t}$$ we have $
%
\begin{remark*}
The point $h=h_{0}, \beta=\beta_{0}:=(3/2)^{3}$ is the critical point and
the curve
%
%e12 ###
\begin{equation}\label{def:gc}
\gamma(t)= \biggl( -\log t - \frac{1+t}{2t}, \frac
{(1+t)^{3}}{2t} \biggr)
\end{equation}
for $t>0$ is the phase transition curve. It corresponds to $\psi
(u^{*})=0$ and $ 2\psi(u^{*})\psi'(u^{*})= 1$. In fact, at the critical
point $(h_{0},\beta_{0})$ the function $\psi(u)=\varphi(u)^{2}-u$
has a
unique root of order three at $u^{*}=4/9$, that is, $\psi(u^{*})=\psi
'(u^{*})=\psi''(u^{*})=0$ and $\psi'''(u^{*})< 0$. The second part of
Lemma \ref{mean2} shows that all the above conclusions (including the
limiting free energy result) are true for the critical point but with
an error rate of $n^{-1/6}$. Define the ``energy'' function
\[
e(r)=\frac{1}{2}I(r,\varphi(0)) + \frac{1}{2} \log\bigl(1-\varphi(0)\bigr)
-\frac
{\beta r^3}{6}
\]
appearing in of the r.h.s. of equation (\ref{eq:lbd}). The ``high
temperature'' regime corresponds to the case when $e(\cdot)$ has a
unique minima and no local maxima or saddle point. The critical point
corresponds to the case when $e(\cdot)$ has a nonquadratic global
minima. The boundary corresponds to the case when $e(\cdot)$ has a
unique minima and a saddle point. In the ``low temperature'' regime,
$e(\cdot)$ has two local minima. In fact, one can easily check that
there is a one-dimensional curve inside the set $S^{c}$, starting from
the critical point, on which $e(\cdot)$ has two global minima and
outside one global minima. Below, we provide the solution on the
boundary curve. Unfortunately, as of now, we don not have a rigorous
solution in the ``low temperature'' regime.

For $(h,\beta)$ on the phase transition boundary curve (excluding the
critical point), the function $\psi(\cdot)$ has two roots and one of
them, say $v^{*}$, is an inflection point. Let $u^{*}$ be the other
root. Here, we mention that $u^{*}$ is a minima of $e(\cdot)$ while
$v^{*}$ is a saddle point of $e(\cdot)$. On the lower part of the
boundary, which corresponds to $\{\gamma(t)\dvtx t<1/2\}$, the inflection
point $v^{*}=(1+t)^{-2}$ is larger than $u^{*}$, while on the upper
part of the boundary corresponding to $\{\gamma(t)\dvtx t>1/2\}$, the
inflection point $v^{*}=(1+t)^{-2}$ is smaller than $u^{*}$. The
following lemma ``solves'' the model at the boundary point $\gamma(t)$
[see (\ref{def:gc})].
\end{remark*}
\begin{lem}\label{lem:boundary}
Let $\gamma(\cdot),u^{*},v^{*}$ be as above and $(h,\beta)=\gamma(t)$
for some
$t\neq1/2$.
Then, for each $1\le i<j\le n$, we have
%
%e13 ###
\begin{equation}
\E(|L_{ij}-u^{*}|)\le\frac{K(\beta,h)}{n^{1/2}}
\end{equation}
for some constant $K(\beta,h)$ depending on $\beta,h$.
Moreover, we have
\[
\frac{\log Z_n(\beta,h)}{n^2} = -\frac{1}{2}I(\varphi(u^*),\varphi(0))
- \frac{1}{2} \log\bigl(1-\varphi(0)\bigr) +\frac{\beta\varphi
(u^*)^3}{6}+O(n^{-1/2})
\]
and
%
%e14 ###
\begin{eqnarray}
&&\pr\left( \left|T_n(\vY)- \pmatrix{n \cr3} \varphi(u^*)^3 \right|
\le
C(\beta,h)n^{5/2} \right)\nonumber\\[-8pt]\\[-8pt]
&&\qquad = \exp\biggl(-\frac{n^2 I(\varphi(u^*),\varphi
(0))}{2}\bigl(1+O(n^{-1/2})\bigr) \biggr),\nonumber
\end{eqnarray}
where $\vY=((Y_{ij}))_{i<j}$ follows G$(n,\varphi(0))$ and the constant
appearing in $O(\cdot)$ and $C(\beta,h)$ depend only on $\beta,h$.
\end{lem}

In the next subsection, we will briefly discuss about the results for
general subgraph counts that can be proved using similar ideas.

%%%%%%%%%%%%%%%%%%%%%%%%%%%%%%%
%s2.3 ###
\subsection{Example: General subgraph counts}\label{subsec:small}

Let $F=(V(F), E(F))$ be a fixed finite graph on $\mvv_{F} :=|V(F)|$
many vertices with $\mve_{F} :=|E(F)|$ many edges. Without loss of
generality, we will assume that $V(F)=[\mvv_{F}]:=\{1,2,\ldots,\mvv
_{F}\}$. Let $\alpha_{F}=|{\operatorname{Aut}}(F)|$ be the number of graph
automorphism of the graph $F$. Let $N_{n}$ be the number of copies of
$F$, not necessarily induced, in the Erd\H{o}s--R\'enyi random graph
$G(n,p)$ (so the number of $2$-stars in a triangle will be three). We
have the following result about the large deviation rate function for
the random variable $N_{n}$.
\begin{theorem}\label{largedevgen}
Let $N_n$ be the number of copies of $F$ in $G(n,p)$, where
\[
p > p_{0}:=\frac{\mve_{ F}-1}{\mve_{F}-1+\exp({\mve
_{F}}/({\mve
_{F}-1}) )}.
\]
Then for any $r\in(p,1]$,
%
%e15 ###
\begin{equation}
\pr\left(N_n \ge\frac{\mvv_{F}!}{\alpha_{F}}\pmatrix{n \cr\mvv_{F}}
r^{\mve
_{F}} \right) = \exp\biggl(-\frac{n^2
I(r,p)}{2}\bigl(1+O(n^{-1/2})\bigr) \biggr).
\end{equation}
Moreover, even if $p \le p_{0}$, there exist $p',p''$ such that
$p<p'\le p''<1$ and the same result holds for all $r\in(p,p')\cup(p'',1]$.
For all $p$ and $r$ in the above domains, we also have the more precise estimate
\begin{eqnarray*}
&&\pr\biggl( \biggl| N_n - \frac{\mvv_{F}!}{\alpha_{F}}\pmatrix{n \cr\mvv_{F}}
r^{\mve_{F}} \biggr| \le C(p,r)n^{\mvv_{F}-1/2} \biggr)\\
&&\qquad = \exp\biggl(-\frac{n^2 I(r,p)}{2}\bigl(1+O(n^{-1/2})\bigr) \biggr),
\end{eqnarray*}
where $C(p,r)$ is a constant depending on $p$ and $r$.
\end{theorem}

\begin{figure}

\includegraphics{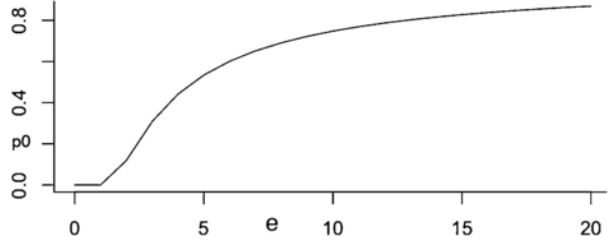}

\caption{The curve $p_{0}$ vs. $\mve_{F}$ where for a graph $F$ with
$\mve_{F}$ many edges our large deviation result holds when $p>p_{0}$.}
\label{fig:p0}
\end{figure}

Note that $p_{0}$ as a function of $\mve_{F}$ is increasing and
converges to $1$ as number of edges goes to infinity (see Figure \ref
{fig:p0}). So there is an obvious gap in the large deviation result,
namely the proof does not work when $r\ge p, p\le p_{0}$ and the gap
becomes larger as the number of edges in $F$ increases. Note that
$p_{0}\to1$ as $\mve_{F}\to\infty$.

%f4 ###

The proof of Theorem \ref{largedevgen} uses the same arguments that
were used in the triangle case. Here, the tilted measure leads to an
exponential random graph model where the Hamiltonian depends on number
of copies of $F$ in the random graph. Let $\beta\ge0, h\in\dR$ be two
fixed numbers. As before, we will identify elements of $\Omega:=\{0,1\}
^{n\choose2}$ with undirected graphs on a set of $n$ vertices. For
each $\mathbf{x}\in\Omega$, let $N(\mathbf{x})$ denote the number
of copies of $F$
in the graph defined by $\mathbf{x}$, and let $E(\mathbf{x})=\sum
_{i<j}x_{ij}$ denote
the number of edges.
Let $\vX= (X_{ij})_{1\le i<j\le n}$ be a random element of $\Omega$
following the probability measure proportional to $e^{H(\mathbf{x})}$, where
$H$ is the Hamiltonian
\[
H(\mathbf{x})=\frac{\beta}{(n-2)_{\mvv_{F}-2}}N(\mathbf{x})+
hE(\mathbf{x}),
\]
where $(n)_{m}=\frac{n!}{(n-m)!}$. Recall that $\mvv_{F}$ is the number
of vertices in the graph $F$. The scaling was done to make the two
summands comparable. Also we used $(n-2)_{\mvv_{F}-2}$ instead of
$n^{\mvv_{F}}$ to make calculations simpler. Let $Z_{n}(\beta,h)$ be the
partition function. Note that $N(\mathbf{x})$ can be written as
%
%e16 ###
\begin{equation}
N(\mathbf{x}) = \frac{1}{\alpha_{F}}\mathop{\sum_{1\le
t_{1},t_{2},\ldots,t_{\mvv
_{F}}\le n,}}_{t_{i}\neq t_{j}\ \mathrm{for}\
i\neq j} \prod_{(i,j)\in E(F)}x_{t_{i}t_{j}}.
\end{equation}
For\vspace*{1pt} $\mathbf{x}\in\Omega, 1\le i<j\le n$, define $\mathbf
{x}_{(i,j)}^{1}$ as the element
of $\Omega$ which is same as $\mathbf{x}$ in every coordinate except
for the
$(i,j)$th coordinate where the value is $1$. Similarly, define $\mathbf{x}
_{(i,j)}^{0}$. For $i<j$, define the random variable
\[
L_{ij}:= \frac{N(\vX_{(i,j)}^{1})-N(\vX_{(i,j)}^{0})}{(n-2)_{\mvv_{F}-2}}.
\]
The main idea is as in the triangle case. We show that $L_{ij}$'s
satisfy a system of ``mean-field equations'' similar to (\ref
{meanfieldeqs}) which has a unique solution under the condition of
Theorem \ref{pressuregen}. In fact, we will show that $L_{ij}$\mbox
{ ``$\approx$'' }$u^{*}$ for all $i<j$ and $E(\vX)$\mbox
{ ``$\approx$'' }${n\choose2}\varphi(u^{*})$ under the condition of Theorem
\ref{pressuregen}. Now note that we can write the Hamiltonian as
\[
H(\vX)= \frac{\beta}{\mve_{F}}\sum_{i<j}X_{ij}L_{ij}+h\sum_{i<j}X_{ij},
\]
which is approximately equal to $ h^{*} E(\vX)$ where $h^{*}=h+\beta
u^{*}/\mve_{F}$. Now the remaining is a calculus exercise.

So the first step in proving the large deviation bound is the following
theorem, which gives the limiting free energy in the ``high
temperature'' regime. Note the similarity with the triangle case.
\begin{theorem}\label{pressuregen}
Suppose we have $\beta\ge0$, $h\in\dR$, and $Z_n(\beta,h)$ defined
as above. Define a function $\varphi\dvtx[0,1]\to\dR$ as
\[
\varphi(x) = \frac{e^{\beta x + h}}{1+e^{\beta x + h}}.
\]
Suppose $\beta$ and $h$ are such that the equation $\alpha_{F} u =
2\mve
_{F}\varphi(u)^{\mve_{F}-1}$ has a unique solution $u^*$ in $[0,1]$ and
$2\mve_{F}(\mve_{F}-1)\varphi(u^{*})^{\mve_{F}-2}\varphi
'(u^{*})<\alpha
_{F}$. Then
\[
\lim_{n\to\infty}\frac{\log Z_n(\beta,h)}{n^2} = -\frac
{1}{2}I(\varphi
(u^*),\varphi(0)) - \frac{1}{2} \log\bigl(1-\varphi(0)\bigr) +\frac{\beta
\varphi
(u^*)^{\mve_{F}}}{\alpha_{F}},
\]
where $I(\cdot,\cdot)$ is the function defined in (\ref{kl}). Moreover,
there exists a constant $K(\beta,h)$ that depends only on $\beta$ and
$h$ (and not on $n$) such that difference between $n^{-2}\log Z_n(\beta
,h)$ and the limit is bounded by $K(\beta,h)n^{-1/2}$ for all $n$.
\end{theorem}

Here also we can identify the region where the conditions in Theorem
\ref{pressuregen} hold. Let
%
%e17 ###
\begin{equation}\label{eq:h0}
h_{0}=\log(\mve_{F}-1)-\frac{\mve_{F}}{\mve_{F}-1}.
\end{equation}
For $h<h_{0}$, there exist exactly two solutions
$0<a_{*}=a_{*}(h)<1/2<a^{*}=a^{*}(h)<\infty$ of the equation
\[
\log x + \frac{1+x}{(\mve_{F}-1)x}+h=0.
\]
Define $a_{*}(h)=a^{*}(h)=1/(\mve_{F}-1)$ for $h=h_{0}$ and
%
%e18 ###
\begin{equation}\label{eq:beta}
\beta_{*}(h)= \frac{\alpha_{F}(1+a_{*})^{\mve_{F}}}{2\mve_{F}(\mve
_{F}-1)a_{*}} \quad\mbox{and}\quad \beta^{*}(h)= \frac{\alpha
_{F}(1+a^{*})^{\mve
_{F}}}{2\mve_{F}(\mve_{F}-1)a^{*}}
\end{equation}
for $h\le h_{0}$.
\begin{lem}\label{lem:pcondgen}
Let $S$ be the set of pairs $(h,\beta)$ for which the function
\[
\psi(u):=2\mve_{F}\varphi(u)^{\mve_{F}-1}-\alpha_{F} u
\]
has a unique root $u^{*}$ in $[0,1]$ and $2\mve_{F}(\mve
_{F}-1)\varphi
(u^{*})^{\mve_{F}-2}\varphi'(u^{*})<\alpha_{F}$ where
$
\varphi(u):= e^{\beta u +h}/(1+e^{\beta u +h}).
$
Then we have
\[
S^{c}=\{(h,\beta)\dvtx h\le h_{0} \mbox{ and } \beta_{*}(h) \le\beta
\le\beta^{*}(h)\},
\]
where $h_{0}, \beta^{*},\beta_{*}$ are as given in (\ref{eq:h0}),
(\ref{eq:beta}).
In particular, $(h,\beta)\in S$ if
\[
\beta\le\frac{\alpha_{F}\mve_{F}^{\mve_{F}-1}}{2(\mve
_{F}-1)^{\mve_{F}}}
\quad\mbox{or}\quad h>h_{0}.
\]
\end{lem}

In fact, Lemma \ref{lem:pcondgen} identifies the critical point and the
phase transition curve where the model goes from ordered phase to a
disordered phase. But the results above does not say what happens at
the boundary or in the low temperature regime. However, note that the
mean-field equations hold for all values of $\beta$ and $h$.

%%%%%%%%%%%%%%%%%%%%%%%%%%%%%%%%%%%%%
%s2.4 ###
\subsection{Example: Ising model on $\bZ^{d}$}
%%%%%%%%%%%%%%%%%%%%%%%%%%%%%%%%%%%%%

Fix any $\beta\ge0, h\in\dR$ and an integer $d\ge1$. Also fix
$n\ge2$.
Let $\bB= \{1,2,\ldots,n+1\}^{d}$ be a hypercube with $(n+1)^{d}$ many
points in the $d$-dimensional hypercube lattice $\bZ^{d}$. Let $\Omega
$ be
the graph obtained from $\bB$ by identifying the opposite boundary
points, that is, for $x=(x_{1},x_{2},\ldots,x_{d}),
y=(y_{1},y_{2},\ldots,y_{d})\in\bB$ we have $x$ is identified with $y$
if $x_{i}-y_{i}\in\{-n,0,n\}$ for all $i$. This identification is known
in the literature as periodic boundary condition. Note that $\Omega$ is
the $d$-dimensional lattice torus with linear size $n$. We will write
$x\sim y$ for $x,y\in\Omega$ if $x,y$ are nearest neighbors in
$\Omega$.
Also, let us denote by $N_{x}$ the set of nearest neighbors of $x$ in
$\Omega$, that is, $N_{x}=\{y\in\Omega\dvtx y\sim x\}$.

Now, consider the Gibbs measure on $\{+1,-1\}^{\Omega}$ given by the
following Hamiltonian
\[
H(\mvgs):= \beta\sum_{x\sim y, x,y\in\Omega} \sigma_{x}\sigma
_{y} + h\sum
_{x\in\Omega
} \sigma_{x},
\]
where $\mvgs=(\sigma_x)_{x\in\Omega}$ is a typical element of $\{
+1,-1\}
^{\Omega
}$. So the probability of a configuration $\mvgs\in\{+1,-1\}^{\Omega
}$ is
%
%e19 ###
\begin{equation} \label{eq:gibbs}
\mu_{\beta,h}(\{\mvgs\})
:= Z_{\beta,h}^{-1}\exp( H(\mvgs) )
=Z_{\beta,h}^{-1}\exp\biggl(\beta\sum_{x\sim y, x,y\in\Omega} \sigma
_{x}\sigma_{y} +
h\sum_{x\in\Omega} \sigma_{x}\biggr),\hspace*{-28pt}
\end{equation}
where $Z_{\beta,h}=\sum_{\mvgs\in\{+1,-1\}^{\Omega}}e^{H(\mvgs)}$
is the
normalizing constant. Here $\sigma_{x}$ is the spin of the magnetic
particle at position $x$ in the discrete torus $\Omega$. This is the
famous Ising model of ferromagnetism on the box $\bB$ with periodic
boundary condition at inverse temperature $\beta$ and external field $h$.

The one-dimensional Ising model is probably the first statistical
model of ferromagnetism to be proposed or analyzed \cite{ising25}. The
model exhibits no phase transition in one dimension. But for dimensions
two and above the Ising ferromagnet undergoes a transition from an
ordered to a disordered phase as $\beta$ crosses a critical value. The
two-dimensional Ising model with no external field was first solved by
Lars Onsager in a ground breaking paper \cite{onsager44}, who also
calculated the critical $\beta$ as $\beta_{c}=\sinh^{-1}(1)$. For
dimensions three and above the model is yet to be solved, and indeed,
very few rigorous results are known.

In this subsection, we present some concentration inequalities for the
Ising model that hold for all values of $\beta$. These
``temperature-free''
relations are analogous to the mean field equations that we obtained
for subgraph counts earlier.

The magnetization of the system, as a function of the configuration
$\mvgs$, is defined as $m(\mvgs):= \frac{1}{|\Omega|}\sum_{x\in
\Omega}
\sigma_{x}$.
%
%Define the average nearest neighbor interaction between the spins in a
%spin configuration $\mvgs$ as
%r(\mvgs)=\frac{1}{d|\Omega|} \sum_{\substack{x\sim y\\x,y\in\Omega}}
%Note that $|\Omega|=n^{d}$ and the Hamiltonian of the spin system is
%simply $n^{d}(\beta d r(\mvgs) + hm(\mvgs))$.
%
For each integer $k\in\{1,2,\ldots,2d\}$, define a degree $k$
polynomial function $r_{k}(\mvgs)$ of a spin configuration $\mvgs$ as follows:
%
%e20 ###
\begin{equation}\label{def:rk}
r_{k}(\mvgs):= \biggl(\pmatrix{2d\cr k}|\Omega| \biggr)^{-1}\sum_{x\in\Omega
}\sum
_{S\subseteq N_{x},|S|=k}\sigma_{S},
\end{equation}
where $\sigma_{S}=\prod_{x\in S}\sigma_{x}$ for any $S\subseteq
\Omega$. In
particular $r_{k}(\mvgs)$ is the average of the product of spins of all
possible $k$ out of $2d$ neighbors.
Note that $r_{1}(\mvgs)\equiv m(\mvgs)$.
% and for $k\ge3$ we have
%m_{k}(\mvgs) = ({2d\choose k}|\Omega| )^{-1}\sum_{
%for
%all } 1\le i\le j\le k}} \prod_{i=1}^{k}\sigma_{y_{i}}.
We will show that when $h = 0$ and $n$ is large, $m(\mvgs)$ and
$r_{k}(\mvgs)$'s satisfy the following ``mean-field relation'' with
high probability under the Gibbs measure:
%
%e21 ###
\begin{equation}\label{eq:mfrel}
\bigl(1-\theta_{0}(\beta)\bigr)m(\mvgs) \approx\sum_{k=1}^{d-1}\theta
_{k}(\beta)
r_{2k+1}(\mvgs).
\end{equation}
These relations hold for all values of $\beta\ge0$. Here, $\theta_{k}$'s
are explicit rational functions of $\tanh(2\beta)$ for $k=0,1,\ldots
,d-1$, defined in (\ref{def:thetak}) below. [Later, we will
prove in Proposition \ref{pr:isingh} that an external magnetic field
$h$ will add an extra linear term in the above relation (\ref
{eq:mfrel}).] The following proposition makes this notion precise in
terms of finite sample tail bound. It is a simple consequence of
Theorem~\ref{thm:conc}.
\begin{theorem}\label{pr:ising}
Suppose $\mvgs$ is drawn from the Gibbs measure $\mu_{\beta,0}$. Then,
for any $\beta\ge0, n\ge1$ and $ t \ge0$ we have
\[
\pr\Biggl( \sqrt{|\Omega|} \Biggl|\bigl(1-\theta_{0}(\beta)\bigr)m(\mvgs) - \sum
_{k=1}^{d-1}\theta_{k}(\beta) r_{2k+1}(\mvgs) \Biggr| \ge t \Biggr)
\le2\exp\biggl( - \frac{t^{2} }{4 b(\beta)} \biggr),
\]
where\vspace*{1pt} $m(\mvgs):= \frac{1}{|\Omega|}\sum_{x\in\Omega} \sigma
_{x}$ is the
magnetization, $r_{k}(\mvgs)$ is as given in (\ref{def:rk}) and for
$k=0,1,\ldots,d-1$
%
%e22 ###
\begin{eqnarray}\label{def:thetak}
\theta_{k}(\beta) &=& \frac{1}{4^{d}}\pmatrix{2d\cr2k+1}\sum
_{\mvgs
\in
\{
-1,+1\}^{2d}} \tanh\Biggl(\beta\sum_{i=1}^{2d}\sigma_{i} \Biggr)\prod
_{j=1}^{2k+1}\sigma_{j}\quad\mbox{and} \nonumber\\[-8pt]\\[-8pt]
b(\beta) &=& |1-\theta_{0}(\beta)| +
\sum_{k=1}^{d-1}(2k+1)|\theta_{k}(\beta)|.\nonumber
\end{eqnarray}
Moreover, we can explicitly write down $\theta_{0}(\beta)$ as
\[
\theta_{0}(\beta) = \frac{1}{4^{d-1}}\sum_{k=1}^{d}k\pmatrix{2d\cr d+k}
\tanh
(2k\beta)
\]
and for $d\ge2$ there exists $\beta_{1}\in(0,\infty)$, depending on
$d$, such that $1-\theta_{0}(\beta)>0$ for $\beta<\beta_{1}$ and
$1-\theta
_{0}(\beta)<0$ for $\beta>\beta_{1}$.
\end{theorem}

%Note that, for dimension two and above, the above result gives a
%transition point for $\beta$ where the magnetization changes sign
%w.r.t.
%the other $r_{k}$'s. from the graphs it seems that $(-1)^{k}\theta_{k}(
%but
%we do not have an analytic proof of that.
%For dimension two the cutoff point is strictly smaller than the
%critical $\beta_{c}$. But in general, we still have no idea whether
%this
%transition has any physical interpretation or it is only a
%computational behavior.
Here, we may remark that for any fixed $k$, $\theta_{k}(\beta/2d)$
converges to the coefficient of $x^{2k+1}$ in the power series
expansion of $\tanh(\beta x)$ and $2d\beta_{1}(d)\downarrow1$ as
$d\to
\infty$. For small values of $d$, we can explicitly calculate the
$\theta_{k}$'s.
%Note that we always have
%$
%$
%In fact we can explicitly write down $\theta_{d-1}(\beta)$ as
%d+k}(-1)^{d+k} \tanh(2k\beta).
%Using the above relations we have the follwing.
For instance, in $d=2$,
\[
\theta_{0}(\beta)=\tfrac{1}{2} \bigl( \tanh(4\beta) + 2\tanh(2\beta
) \bigr),\qquad
\theta_{1}(\beta)=\tfrac{1}{2} \bigl( \tanh(4\beta) - 2\tanh(2\beta
) \bigr).
\]
For $d=3$,
\begin{eqnarray*}
\theta_{0}(\beta)&=&\tfrac{3}{16} \bigl( \tanh(6\beta) + 4\tanh(4\beta
)+5\tanh
(2\beta) \bigr),\\
\theta_{1}(\beta)&=&\tfrac{10}{16} \bigl( \tanh(6\beta) - 3\tanh(2\beta
)
\bigr),\\
\theta_{2}(\beta)&=&\tfrac{3}{16} \bigl( \tanh(6\beta) - 4\tanh(4\beta
)+5\tanh
(2\beta) \bigr).
\end{eqnarray*}
For $d=4$,
\begin{eqnarray*}
\theta_{0}(\beta)&=&\tfrac{1}{16} \bigl(\tanh(8\beta) + 6\tanh(6\beta) +
14\tanh
(4\beta)+14\tanh(2\beta) \bigr),\\
\theta_{1}(\beta)&=&\tfrac{7}{16} \bigl(\tanh(8\beta) + 2\tanh(6\beta) -
2\tanh
(4\beta) - 6\tanh(2\beta) \bigr),\\
\theta_{2}(\beta)&=&\tfrac{7}{16} \bigl(\tanh(8\beta) - 2\tanh(6\beta) -
2\tanh
(4\beta)+6\tanh(2\beta) \bigr),\\
\theta_{3}(\beta)&=&\tfrac{1}{16} \bigl(\tanh(8\beta) - 6\tanh(6\beta) +
14\tanh
(4\beta)-14\tanh(2\beta) \bigr).
\end{eqnarray*}
%
%Note that
%where $u=\tanh(2\beta)$.

\begin{cor}
For the Ising model on $\Omega$ at inverse temperature $\beta$ with no
external magnetic field for all $ t\ge0$ we have:
\begin{longlist}
\item if $d=1$,
\[
\pr\bigl( |m(\mvgs)|\ge t\bigr) \le2\exp\bigl( - \tfrac{1}{4}|\Omega|\bigl(1-\tanh
(2\beta
)\bigr)t^{2} \bigr);
\]
\item if $d=2$,
\[
\pr\bigl( |[(1-u)^{2}-u^{3}]m(\mvgs) + u^{3}r_{3}(\mvgs)|\ge t\bigr)\le2\exp
\biggl( - \frac{|\Omega|t^{2}}{32} \biggr),
\]
where $u=\tanh(2\beta)$ and $r_{3}(\mvgs)= \frac{1}{4|\Omega|}\sum^{*}
\sigma
_{x}\sigma_{y}\sigma_{z}$ where the sum $\sum^{*}$ is over all
$x,y,z\in
\Omega$
such that $|x-y|=2,|z-y|=2,|x-z|=2$;
\item if $d=3$,
\[
\pr\bigl( |g(u)m(\mvgs) + 5u^{3}(1+u^{2})r_{3}(\mvgs) -3u^{5}r_{5}(\mvgs
)|\ge t\bigr)
\le2\exp( - c|\Omega|t^{2} ),
\]
where $c$ is an absolute constant, $g(u) = 1-3u
+4u^{2}-9u^{3}+3u^{4}-3u^{5}$, $u=\tanh(2\beta)$ and $r_{3}, r_{5}$ are
as defined in (\ref{def:rk}).
\end{longlist}
\end{cor}

%Though the results obtained for Ising model on hypercube lattice are
%easy consequences of Theorem \ref{thm:conc}, we haven't encountered
%any explicit result like the above in the literature. Note that, the
%relations are valid for all values of $\beta$. In fact, a similar
%calculation can be used to prove the same result for the hypercube $
%allows the same proof to hold is that $|\partial\bB|=o(|\bB|)$.

Although we do not yet know the significance of the above relations, it
seems somewhat striking that they are not affected by phase
transitions. The exponential tail bounds show that many such relations
can hold simultaneously.
For completeness, we state below the corresponding result for nonzero
external field.

\begin{prop}\label{pr:isingh}
Suppose $\mvgs$ is drawn from the Gibbs measure $\mu_{\beta,h}$. Let
$r_{k}(\mvgs)$, $\theta_{k}(\beta)$, $b(\beta)$ be as in proposition
(\ref{pr:ising}). Then, for any $\beta\ge0, h\in\dR, n\ge1$ and $ t
\ge0$
we have
%
%e23 ###
\begin{equation}
\pr\bigl( \bigl|\bigl(1-\theta_{0}(\beta)\bigr)m(\mvgs) - g(\mvgs) \bigr|
\ge
t \bigr)\le2\exp\biggl( - \frac{|\Omega|t^{2} }{4 b(\beta)(1+\tanh|h|)}
\biggr),
\end{equation}
where
\[
g(\mvgs):= \sum_{k=1}^{d-1}\theta_{k}(\beta) r_{2k+1}(\mvgs)
+\tanh(h)
\Biggl(1-\sum_{k=0}^{d-1}\theta_{k}(\beta) s_{2k+1}(\mvgs) \Biggr)
\]
and
\[
s_{k}(\mvgs):= \biggl(\pmatrix{2d\cr k}|\Omega| \biggr)^{-1}\sum_{x\in\Omega
}\sum
_{S\subseteq N_{x}, |S|=k}\sigma_{S\cup\{x\}}
\]
is the average of products of spins over all $k$-stars for
$k=1,2,\ldots
,2d$ and $\Omega$ is the discrete torus in $\bZ^{d}$ with $n^{d}$
many points.\vadjust{\goodbreak}
\end{prop}
\eject
%
%%%%%%%%%%%%%%%%%%%%%%%%%%%%%%%%%%%%%

%%%%%%%%%%%%%%%%%%%%%%%%%%%%%%%%%%%%%
%s3 ###
\section{Proofs}\label{sec:proof}

%%%%%%%%%%%%%%%%%%%%%%%%%%%%%%%%%%%%%%%%%%%%
%s3.1 ###
\subsection{\texorpdfstring{Proof of Proposition
\protect\ref{pr:crit}}{Proof of Proposition 4}}

Instead of proving Theorem \ref{thm:nong1} first, let us see how it is
applied to prove the result for the Curie--Weiss model at critical
temperature. The proof is simply an elaboration of the sketch given at
the end of Section \ref{subsec:curie}.

Suppose $\mvgs$ is drawn from the Curie--Weiss model at critical
temperature. We construct $\mvgs'$ by taking one step in the heat-bath
Glauber dynamics: a coordinate $I$ is chosen uniformly at random, and
$\sigma_{I}$ is replace by $\sigma'_{I}$ drawn from the conditional
distribution of the $I$th coordinate given $\{\sigma_{j}\dvtx j\neq I\}$. Let
\[
F(\mvgs,\mvgs') := \sum_{i=1}^{n}(\sigma_{i}-\sigma'_{i})=\sigma
_{I}-\sigma'_{I}.
\]
For each $i=1,2,\ldots,n$, define
$
m_{i}=m_{i}(\mvgs) = n^{-1}\sum_{j\neq i} \sigma_{j}.
$
An easy computaion gives that $\E(\sigma_{i}| \{\sigma_{j}, j\neq i\}
) =
\tanh( m_{i})$ for all $i$ and so we have
\[
f(\mvgs):= \E(F(\mvgs,\mvgs')|\mvgs) = m - \frac{1}{n} \sum_{i=1}^{n}
\tanh(m_{i})=\frac{m}{n} + \frac{1}{n} \sum_{i=1}^{n} g(m_{i}),
\]
where $g(x):=x-\tanh(x)$.
By definition, $m_{i}(\mvgs)-m(\mvgs) = \sigma_{i}/n$ and
$m_{i}(\mvgs
')-m(\mvgs) = (\sigma_{i}+\sigma_{I}-\sigma'_{I})/n$ for all $i$.
Hence, using
Taylor's expansion up to first degree and noting that $|g'(x)|=\tanh
^{2}(x)\le x^{2}$ we have
\begin{eqnarray*}
|f(\mvgs)-f(\mvgs')| &\le& \frac{2}{n}|g'(m(\mvgs))|+ \frac{2+5\max
_{|x|\le1}|g''(x)|}{n^{2}}\\
&\le& \frac{2}{n}m(\mvgs)^{2} +\frac{6}{n^{2}}.
\end{eqnarray*}
Clearly, $|F(\mvgs,\mvgs')|\le2$. Thus, we have
\[
\Delta(\mvgs) := \frac{1}{2}\E[ |f(\mvgs)-f(\mvgs')|\cdot
|F(\mvgs
,\mvgs')| \mid\mvgs]
\le\frac{2}{n}m(\mvgs)^{2}+\frac{6}{n^{2}}.
\]
Now it is easy to verify that $|x|^{3} \le5|x-\tanh x|$ for all
$|x|\le1$. Note that this is the place where we need $\beta=1$. For
$\beta
\neq1$, the linear term dominates in $m-\tanh(\beta m)$. Hence, it
follows that
\[
m(\mvgs)^{2} \le5^{2/3}|m (\mvgs)-\tanh m(\mvgs) |^{2/3}\le
3|f(\mvgs
)|^{2/3} +{3}{n^{-2/3}},
\]
where in the last line we used the fact that $|f(\mvgs) - (m - \tanh
m)|\le1/n$ and $5^{2/3}< 3$. Thus,
\[
\Delta(\mvgs) \le\frac{6}{n}|f(\mvgs)|^{2/3}+\frac{12}{n^{5/3}}
\]
and using Corollary \ref{cor:nong} with $\alpha=2/3, B= 6/n$ and
$C=12/n^{5/3}$ we have
\[
\pr(|m-\tanh m| \ge t+n^{-1}) \le\pr\bigl(|f(\mvgs)| \ge t\bigr) \le2e^{- cnt^{4/3}}
\]
for all $t\ge0$ for some constant $c>0$.
This clearly implies that
\[
\pr(|m|\ge t) \le\pr(|m-\tanh m|\ge t^{3}/5) \le2e^{-cnt^{4}}
\]
for all $t\ge0$ and for some absolute constant $c>0$. Thus, we are done.

%%%%%%%%%%%%%%%%%%%%%%%%%%%%%%%%%%%%
%s3.2 ###
\subsection{\texorpdfstring{Proof of Proposition
\protect\ref{pr:critgen}}{Proof of Proposition 6}}

The proof is along the lines of proof of Proposition \ref{pr:crit}.
Suppose $\vX$ is drawn from the distribution $\nu_{n}$. We construct
$\vX'$ as follows: a coordinate $I$ is chosen uniformly at random, and
$X_{I}$ is replace by $X'_{I}$ drawn from the conditional distribution
of the $I$th coordinate given $\{X_{j}\dvtx j\neq I\}$. Let
\[
F(\vX,\vX') := \sum_{i=1}^{n}(X_{i}-X'_{i})=X_{I}-X'_{I}.
\]
For each $i=1,2,\ldots,n$, define
$
m_{i}(\vX) = n^{-1}\sum_{j\neq i} X_{j}.
$
An easy computaion gives that $\E(X_{i}| \{X_{j}, j\neq i\}) =
g(m_{i})$ for all $i=1,2,\ldots,n$ where $g(s)=\frac{d}{ds}(\log\int
\exp(x^{2}/2n+sx) \,d\rho(x))$ for $s\in\dR$. So we have
\[
f(\vX):= \E(F(\vX,\vX')|\vX) = m(\vX) - \frac{1}{n} \sum_{i=1}^{n}
g(m_{i}(\vX)).
\]
Define the function
%
%e24 ###
\begin{equation}\label{def:h}
h(s)=\frac{s^{2}}{2}-\log\int\exp(sx) \,d\rho(x) \qquad\mbox{for }s\in
\dR.
\end{equation}
Clearly, $h$ is an even function. Recall that $k$ is an integer such
that $h^{(i)}(0)=0$ for $0\le i<2k$ and $h^{(2k)}(0)\neq0$. We have
$k\ge2$ since $h''(0)=1-\int x^{2} \,d\rho(x)=0$.

Now using the fact that $\rho([-L,L])=1$ it is easy to see that
$|f(\vX
)-h'(m(\vX))|\le c/n$ for some constant $c$ depending on $L$ only. In
the subsequent calculations, $c$ will always denote a constant
depending only on $L$ that may vary from line to line. Similarly, we have
\begin{eqnarray*}
|f(\vX)-f(\vX')|&\le&\frac{|X_{I}-X'_{I}|}{n} \biggl( |1-g'(m(\vX
))|+\frac
{c(1+\sup_{|x|\le L}|g''(x)|)}{n} \biggr)\\
&\le&\frac{2L}{n}|h''(m(\vX))| + \frac{c}{n^{2}}.
\end{eqnarray*}
Note that\vspace*{1pt} $|h''(s)|\le cs^{2k-2}$ for some constant $c$ for all $s\ge
0$. This follows since $\lim_{s\to0}h''(s)/s^{2k-2}$ exists\vspace*{1pt} and
$h''(\cdot)$ is a bounded function. Also $\lim_{s\to
0}|h'(s)|/\break|s|^{2k-1}=|h^{(2k)}(0)|\neq0$ and $|h'(s)|> 0$ for $s> 0$.
So we have $|h'(s)|\ge c|s|^{2k-1}$ for some constant $c>0$ and all
$|s|\le L$. From the above results, we deduce that
\begin{eqnarray*}
|f(\vX)-f(\vX')|&\le&\frac{c}{n}|(m(\vX))|^{2k-2} + \frac{c}{n^{2}}
\le\frac{c}{n}|h'(m(\vX))|^{({2k-2})/({2k-1})} + \frac{c}{n^{2}}\\
&\le&\frac{c}{n}|f(\vX)|^{({2k-2})/({2k-1})} + \frac{c}{n^{2-{1}/{(2k-1)}}}.
\end{eqnarray*}
Now the rest of the proof follows exactly as for the classical
Curie--Weiss model.

%%%%%%%%%%%%%%%%%%%%%%%%%%%%%%%%%%%%%

%s3.3 ###
\subsection{\texorpdfstring{Proof of Theorem \protect\ref{largedev}}{Proof of Theorem 7}}

First, let us state and prove a simple technical lemma.
\begin{lem}\label{techlmm}
Let $x_1,\ldots,x_k,y_1,\ldots,y_k$ be real numbers. Then
\[
\max_{1\le i\le n} \biggl|\frac{e^{x_i}}{\sum_{j=1}^k e^{x_j}} -
\frac
{e^{y_i}}{\sum_{j=1}^k e^{y_j}} \biggr|\le2\max_{1\le i\le n} |x_i - y_i|
\]
and
\[
\Biggl|\log\sum_{i=1}^k e^{x_i} - \log\sum_{i=1}^k e^{y_i}
\Biggr|\le
\max_{1\le i\le k} |x_i - y_i|.
\]
\end{lem}
\begin{pf}
Fix $1\le i\le k$. For $t\in[0,1]$, let
\[
h(t) = \frac{e^{tx_i + (1-t)y_i}}{\sum_{j=1}^k e^{tx_j + (1-t)y_j}}.
\]
Then
\[
h'(t) = \biggl[(x_i -y_i) - \frac{\sum_{j=1}^k (x_j-y_j) e^{tx_j +
(1-t)y_j}}{\sum_{j=1}^ke^{tx_j + (1-t)y_j}} \biggr] h(t).
\]
This shows that $|h'(t)| \le2\max_i |x_i-y_i|$ for all $t\in[0,1]$
and completes the proof of the first assertion. The second inequality
is proved similarly.
\end{pf}
\begin{pf*}{Proof of Lemma \ref{lem:mf}}
Fix two numbers $1\le i < j\le n$. Given a configuration~$\vX$,
construct another configuration $\vX'$ as follows. Choose a point $k
\in\{1,\ldots,n\}\setminus\{i,j\}$ uniformly at random, and replace
the pair $(X_{ik}, X_{jk})$ with $(X'_{ik}, X'_{jk})$ drawn from the
conditional distribution given the rest of the edges. Let $L_{ij}'$ be
the revised value of $L_{ij}$. From the form of the Hamiltonian, it is
now easy to read off that for $x,y\in\{0,1\}$,
\begin{eqnarray*}
\hspace*{-3pt}&&\pr(X'_{ik} = x, X'_{jk} = y \mid\vX) \\
\hspace*{-3pt}&&\qquad\propto\exp\biggl( \beta x L_{ik} + \beta y L_{jk} + hx + hy -
\frac
{\beta}{n} x X_{ij} X_{jk} - \frac{\beta}{n} y X_{ij} X_{ik} + \frac
{\beta}{n} xyX_{ij} \biggr).
\end{eqnarray*}
An application of Lemma \ref{techlmm} shows that the terms having
$\beta
/n$ as coefficient can be ``ignored'' in the sense that for each $x,y\in
\{0,1\}$,
\[
\biggl|\pr(X'_{ik} = x, X'_{jk} = y \mid\vX) - \frac{e^{\beta x L_{ik}
+ \beta y L_{jk} + hx + hy}}{(1+e^{\beta L_{ik} + h}) (1+e^{\beta
L_{jk} + h})} \biggr|\le\frac{2\beta}{n}.
\]
In particular,
%
%e25 ###
\begin{equation}\label{diff}
|\E(X'_{ik} X'_{jk} \mid\vX) - \varphi(L_{ik})\varphi(L_{jk})| \le
\frac{2\beta}{n}.
\end{equation}
Now,
%
%e26 ###
\begin{eqnarray}\label{fx}
\E(L_{ij} - L_{ij}'\mid\vX) &=& \frac{1}{n(n-2)} \sum_{k\notin\{
i,j\}
} \bigl(X_{ik}X_{jk} - \E(X'_{ik}X'_{jk}\mid\vX)\bigr)\nonumber\\[-8pt]\\[-8pt]
&=& \frac{1}{n-2} L_{ij} - \frac{1}{n(n-2)}\sum_{k\notin\{i,j\}} \E
(X'_{ik} X'_{jk}\mid
\vX).\nonumber
\end{eqnarray}
Let $F(\vX,\vX') = (n-2) (L_{ij} - L'_{ij})$ and $f(\vX) = \E(F(\vX
,\vX
')\mid\vX)$. Let
\[
g(\vX) = L_{ij} - \frac{1}{n}\sum_{k\notin\{i,j\}}\varphi
(L_{ik})\varphi(L_{jk}).
\]
From (\ref{diff}) and (\ref{fx}), it follows that
%
%e27 ###
\begin{equation}\label{eq:fg}
|f(\vX) - g(\vX)| \le\frac{2\beta}{n}.
\end{equation}
Since $X'$ has the same distribution as $X$, the same bound holds for
$|f(X')-g(X')|$ as well. Now clearly, $|F(X,X')| \le1$. Again, $|g(X)
- g(X')| \le2/n$, and therefore
\[
|f(X) - f(X')|\le\frac{4(1+\beta)}{n}.
\]
Combining everything, and applying Theorem \ref{thm:conc} with $B=0$
and $C=2(1+\beta)/n$, we get
\[
\pr\bigl(|f(\vX)|\ge t\bigr)\le2\exp\biggl( - \frac{nt^{2}}{4(1+\beta)} \biggr)
\]
for all $t\ge0$. From (\ref{eq:fg}), it follows that
\[
\pr\bigl(|g(\vX)|\ge t\bigr)\le\pr\bigl(|f(\vX)|\ge t-2\beta/n\bigr)\le
2\exp\biggl( - \frac{nt^{2}}{8(1+\beta)} \biggr)
\]
for all $t\ge8\beta/n$. This completes the proof of the tail bound. The
bound on the mean absolute value is an easy consequence of the tail bound.
\end{pf*}
\begin{pf*}{Proof of Lemma \ref{mean2}}
The proof is in two steps. In the first step, we will get an error
bound of order $n^{-1/2}\sqrt{\log n}$. In the second step, we will
improve it to $n^{-1/2}$.
Define
\[
\Delta= \max_{1\le i<j\le n} \biggl|L_{ij} - \frac{1}{n}\sum_{k\notin\{
i,j\}} \varphi(L_{ik})\varphi(L_{jk}) \biggr|.
\]
By Lemma \ref{lem:mf} and union bound, we have
\[
\pr(
\Delta\ge t ) \le n^{2}\exp\biggl(- \frac{nt^{2}}{8(1+\beta)}
\biggr)
\]
for all $t\ge8\beta/n$. Intuitively, the above equation says that
$\Delta$
is of the order of $\sqrt{\log n/n}$, in fact we have $\E( \Delta
^{2} )=
O(\log n/n)$. Clearly, $\varphi$ is an increasing function. Hence, we have
\[
\varphi(L_{\min})^{2} - \Delta\le L_{\min}\le L_{\max}\le\varphi
(L_{\max
})^{2} + \Delta,
\]
where $L_{\max}=\max_{1\le i<j\le n}L_{ij}$ and $L_{\min}=\min
_{1\le
i<j\le n}L_{ij}$.

% \centering
% % \fbox{\begin{minipage}{3in}
% \includegraphics[width=4in]{3.jpg}
% \caption{Choice of $\eps=0.05,\gd=0.775$ when $\beta=4,h=-0.5$ and
%$u^{*}=0.9225$.}
% \label{fig:ex}

Now assume that there exists a unique solution $u^{*}$ of the equation
$\varphi(u)^{2}=u$ with $2\varphi(u^{*})\varphi'(u^{*})<1$. For ease of
notation, define the function $\psi(u)=\varphi(u)^2-u$. We have
$\psi
(0)>0>\psi(1)$, $u^{*}$ is the unique solution to $\psi(u)=0$ and
$\psi
'(u^{*})<0$. It is easy to see that $\psi'(u)=0$ has at most three
solution [$\psi'(u)=2\beta\varphi(u)^{2}(1-\varphi(u))-1$ is a third
degree polynomial in $\varphi(u)$ and $\varphi$ is a strictly
increasing function].

Hence, there exist positive real numbers $\eps, \delta$ such that
$|\psi
(u)|>\eps$ if $|u-u^{*}|>\delta$. Note that $\psi(u)>0$ if $u<u^{*}$ and
$\psi(u)<0$ is $u>u^{*}$. Decreasing $\eps,\delta$ without loss of
generality, we can assume that
%
%e28 ###
\begin{equation}\label{eq:deriv}
\inf_{0<|u-u^{*}|\le\delta} \biggl[\frac{u-u^{*}}{-\psi(u)}
\biggr]= c >0.
\end{equation}
This is possible because $\psi'(u^{*})<0$. Note that $\psi(L_{\max
})\ge
-\Delta$ and $\psi(L_{\min})\le\Delta$. Thus, we have
\[
u^{*}-\delta\le L_{\min}\le L_{\max}\le u^{*}+\delta,
\]
when $\Delta<\eps$. Using (\ref{eq:deriv}), $u^{*}\le L_{\max}\le
u^{*}+\delta
$ implies that $|L_{\max}-u^{*}|\le c\Delta$ and $u^{*}-\delta\le
L_{\min}\le
u^{*}$ implies that $|L_{\min}-u^{*}|\le c\Delta$.
Thus, when $\Delta< \eps$, we have $|L_{\max}-u^{*}|\le c{\Delta}$ and
$|L_{\min}-u^{*}|\le c{\Delta}$ and in particular, $|L_{ij}-u^{*}|\le
c{\Delta
}$ for all $i<j$. So we can bound the $L^{2}$ distance of $L_{ij}$ from
$u^{*}$ by
\[
\E(L_{ij}-u^{*})^{2}\le c^{2}\E(\Delta^{2}) + \pr(\Delta\ge\eps
)\le
K(\beta
,h)\frac{\log n}{n}
\]
for all $i<j$.

Now let us move to the second step. Recall from (\ref{eq:lcon}) that
%
%e29 ###
\begin{equation}\label{eq:lconcopy}
\E\biggl|L_{ij} - \frac{1}{n}\sum_{k\notin\{i,j\}} \varphi
(L_{ik})\varphi(L_{jk}) \biggr| \le\frac{C(1+\beta)^{1/2}}{n^{1/2}}
\end{equation}
for all $i<j$. Let $D_{ij}=L_{ij}-u^{*}$. Using Taylor's expansion around
$u^{*}$ up to degree one, we have
\begin{eqnarray*}
\varphi(L_{ik})\varphi(L_{jk}) - \varphi(u^{*})^{2}
&=& \varphi(u^{*})\bigl(\varphi(L_{ik}) - \varphi(u^{*})\bigr) +\varphi
(u^{*})\bigl(\varphi(L_{jk}) - \varphi(u^{*})\bigr)\\
&&{} + \bigl(\varphi(L_{ik}) - \varphi(u^{*})\bigr)\bigl(\varphi(L_{jk}) -
\varphi
(u^{*})\bigr) \\
&=& \varphi(u^{*}) \varphi'(u^{*})(D_{ik}+D_{jk}) + R_{ijk},
\end{eqnarray*}
where $\E(|R_{ijk}|)\le C\E(D_{ij}^{2})\le C n^{-1}\log n$ for some
constant $C$ depending only on $\beta,h$.
Thus,
%
%e30 ###
\begin{eqnarray}\label{eq:combl2}
&&\E\biggl|L_{ij} - \frac{1}{n}\sum_{k\notin\{i,j\}} \varphi
(L_{ik})\varphi(L_{jk}) - D_{ij} + \frac{\varphi(u^{*}) \varphi
'(u^{*})}{n}\sum_{k\notin\{i,j\}}(D_{ik}+D_{jk}) \biggr|\hspace*{-28pt}\nonumber\\[-8pt]\\[-8pt]
&&\qquad\le\frac{2u^{*}}{n} +\frac{1}{n}\sum_{k\notin\{i,j\}}\E
|R_{ijk}|\le\frac{C\log
n}{n}.\nonumber
\end{eqnarray}
Here, we used the fact that $u^{*}=\varphi(u^{*})^{2}$. Combining
(\ref{eq:lconcopy}) and (\ref{eq:combl2}), we have
\[
\E\biggl|D_{ij} - \frac{\varphi(u^{*}) \varphi'(u^{*})}{n}\sum
_{k\notin\{i,j\}}(D_{ik}+D_{jk}) \biggr|\le\frac{C}{\sqrt{n}}
\]
for all $i<j$. By symmetry, $\E|D_{ij}|$ is the same for all $i,j$.
Thus, finally we have
\[
\E|L_{ij}-u^{*}|=\E|D_{ij}|\le\frac{1}{1-2\varphi(u^{*}) \varphi
'(u^{*})}\cdot\frac{C}{\sqrt{n}}=\frac{K(\beta,h)}{\sqrt{n}},
\]
where $K(\beta,h)$ is a constant depending on $\beta,h$.

When $\psi(u)=0$ has a unique solution at $u=u^{*}$ with $2\psi
(u^{*})\psi'(u^{*})=1$, which happens at the critical point $\beta
=(3/2)^{3}, h=\log2 -3/2$, instead of (\ref{eq:deriv}) we have
\[
\inf_{0<|u-u^{*}|\le\delta} \biggl[\frac{(u-u^{*})^{3}}{-\psi
(u)} \biggr]=
c >0
\]
since $\psi(u^{*})=\psi'(u^{*})=\psi''(u^{*})=0$ and $\psi
'''(u^{*})<0$. Then using a similar idea as above one can easily show that
\[
\E|L_{ij}-u^{*}|\le K(\beta,h){n^{-1/6}}
\]
for some constant $K$ depending on $\beta,h$. This completes the proof of
the lemma.
%\rightqed
\end{pf*}
\begin{remark*}
The proof becomes lot easier if we have
%
%e31 ###
\begin{equation}\label{ccond}
c:=\varphi(1)\cdot\sup_{0\le x\le1}\frac{|\varphi(x) - \varphi
(u^*)|}{|x-u^*|}<\frac{1}{2}.
\end{equation}
This is because, by the triangle inequality, we have
%
%e32 ###
\begin{eqnarray}\label{maineq}
\sum_{i<j} |L_{ij} - u^*| &\le& \sum_{i<j} \biggl|L_{ij} - \frac
{1}{n}\sum_{k\notin\{i,j\}}\varphi(L_{ik})\varphi(L_{jk})
\biggr|\nonumber\\[-8pt]\\[-8pt]
&&{} + \sum_{i<j} \biggl(\frac{1}{n}\sum_{k\notin\{i,j\}}
|\varphi(L_{ik})\varphi(L_{jk}) - u^*| +
\frac{2u^*}{n} \biggr).\nonumber
\end{eqnarray}
Now recall that condition (\ref{ccond}) says that $\varphi(1)|\varphi
(x)- \varphi(u^*)|\le c|x-u^*|$ for all $x\in[0,1]$. Moreover,
$L_{ij}\in[0,1]$ for all $i,j$, and $u^* = \varphi(u^*)^2$. Thus,
\[
|\varphi(L_{ik})\varphi(L_{jk}) - u^*| \le c|L_{ik} -
u^*| +
c|L_{jk} - u^*|.
\]
Combining everything, we get
\[
\sum_{i<j} |L_{ij} - u^*| \le\frac{\sum_{i<j}|L_{ij} -
{1}/{n}\sum_{k\notin\{i,j\}}\varphi(L_{ik})\varphi(L_{jk})|
+ n
u^*}{1-2c}.
\]
Taking expectation on both sides, and applying Lemma \ref{lem:mf}, we get
\[
\sum_{i<j} \E|L_{ij} - u^*| \le\frac{C(1+\beta) n^{3/2}}{1-2c}.
\]
And this gives the required result. In fact, using basic calculus
results one can easily check that condition (\ref{ccond}) is satisfied
when $h\ge0$ or $\beta\le2$.
\end{remark*}

Now, we will prove that in the exponential random graph model, the
number of edges and number of triangles also satisfy certain
``mean-field'' relations.
\begin{lem}\label{lem:ex}
Recall that $E(\mathbf{x})$ and $T(\mathbf{x})$ denote the number of
edges and number
of triangles in the graph defined by the edge configuration $\mathbf
{x}\in\Omega
$. If $\vX$ is drawn from the Gibbs' measure in Theorem \ref{pressure},
we have the bound
\begin{eqnarray*}
\E\biggl| E(\vX) - \sum_{i<j}\varphi(L_{ij}) \biggr| &\le& C(1+\beta
)^{1/2}n,\\
\E\biggl| \frac{T(\vX)}{n} - \frac{1}{3}\sum_{i<j}L_{ij}\varphi
(L_{ij}) \biggr| &\le& C(1+\beta)^{1/2}n,
\end{eqnarray*}
where and $C$ is a universal constant.
\end{lem}
\begin{pf}
It is not difficult to see that
\[
\E\bigl(X_{ij}\mid(X_{kl})_{ (k,l)\ne(i,j)}\bigr) = \varphi(L_{ij}).
\]
Let us create $\vX'$ by choosing $1\le i<j\le n$ uniformly at random
and replacing $X_{ij}$ with $X_{ij}'$ drawn from the conditional
distribution of $X_{ij}$ given $(X_{kl})_{(k,l)\ne(i,j)}$. Let $F(\vX
,\vX') = {n\choose2}(X_{ij}-X'_{ij})$. Then
\[
f(\vX) = \E(F(\vX,\vX')|\vX) = \sum_{k<l} \bigl(X_{kl} - \varphi
(L_{kl})\bigr) =
E(\vX) - \sum_{k<l}\varphi(L_{kl}).
\]
Now $|F(\vX,\vX')| \le{n \choose2}$ and $|f(\vX) - f(\vX')| \le
1+\beta$. Here we used the fact that $|\varphi'(x)|\le\beta/4$.
Combining the above result and Theorem \ref{thm:conc} with $B=0,
C=\frac
{1}{2}(1+\beta){n\choose2}$, we get the required bound.

Similarly, if we define $F(\vX,\vX') = {n\choose2}(X_{ij}L_{ij} -
X_{ij}'L_{ij})$. Then
\begin{eqnarray*}
f(\vX) &=& \E(F(\vX,\vX')|\vX) = \sum_{k<l} \bigl(X_{kl}L_{kl} -
\varphi
(L_{kl})L_{kl}\bigr)\\
&=& \frac{3}{n}T(\vX) - \sum_{k<l}\varphi(L_{kl})L_{kl}.
\end{eqnarray*}
Again, $|F(\vX,\vX')| \le{n\choose2}$ and $|f(\vX)-f(\vX')| \le
C(1+\beta)$. The bound follows easily as before.
\end{pf}

The following result is an easy corollary of Lemmas \ref{mean2} and
\ref{lem:ex}.
\begin{cor}\label{cor:ex}
Suppose the conditions of Theorem \ref{pressure} are satisfied. Then
we have
\[
\E\biggl\vert E(\vX) - \frac{n^{2}\varphi(u^{*})}{2} \biggr\vert\le Cn^{3/2}
\quad\mbox{and}\quad
\E\biggl\vert\frac{T(\vX)}{n} - \frac{n^{2}\varphi(u^{*})^{3}}{6}
\biggr\vert\le Cn^{3/2},
\]
where $C$ is a constant depending only on $\beta,h$.
\end{cor}
\begin{lem}\label{large1}
Suppose the conditions of Theorem \ref{pressure} are satisfied. Let
$T_n$ be the number of triangles in the Erd\H os--R\'enyi graph $G(n,
\varphi(0))$. Then there is a constant $K(\beta,h)$ depending only on
$\beta$ and $h$ such that for all $n$
\[
\biggl|\frac{\log\pr(|T_n - {n \choose 3}\varphi(u^*)^3| \le
K(\beta
,h) n^{5/2})}{n^2} - \frac{-I(\varphi(u^*),\varphi(0))}{2}
\biggr|\le
\frac{K(\beta,h)}{\sqrt{n}}.
\]
\end{lem}
%
%%%%%%%%%%%%%%%%%%%%%%%%%%%%%%%%%%%%%%%%
%
\begin{pf}
Let $X$ be drawn from the Gibbs' measure in Theorem \ref{pressure} with
parameters $\beta,h$. From Corollary \ref{cor:ex}, we see that there
exists a constant $K(\beta,h)$ such that (for all $n$)
\[
\pr\biggl( \biggl|E(X) - \frac{n^2\varphi(u^*)}{2} \biggr|\le
K(\beta,h)
n^{3/2} \biggr) \ge\frac{3}{4}
\]
and
\[
\pr\biggl( \biggl|\frac{T(X)}{n} - \frac{n^2\varphi(u^*)^3}{6}
\biggr|\le K(\beta,h) n^{3/2} \biggr) \ge\frac{3}{4}.
\]
Now let
\[
A = \biggl\{x\in\{0,1\}^n \dvtx \biggl|\frac{T(x)}{n} - \frac
{n^2\varphi
(u^*)^3}{6} \biggr|\le K(\beta,h) n^{3/2} \biggr\}
\]
and
\[
B = A \cap\biggl\{x\in\{0,1\}^n\dvtx \biggl|E(x) - \frac{n^2\varphi
(u^*)}{2} \biggr|\le K(\beta,h) n^{3/2} \biggr\}.
\]
Now suppose $Y=(Y_{ij})_{1\le i<j\le n}$ is a collection of i.i.d.
random variables satisfying $\pr(Y_{ij} = 1) = 1 - \pr(Y_{ij} = 0) =
\varphi(0)$ and $Z=(Z_{ij})_{1\le i<j\le n}$ is another collection of
i.i.d. random variables with $\pr(Z_{ij}= 1)=1-\pr(Z_{ij}=0) =
\varphi
(u^*)$. Without loss of generality, we can assume that $K(\beta,h)$ was
chosen large enough to ensure that (again, for all $n$) $\pr(Z\in A)
\ge1/2$ and $\pr(Z\in B) \ge1/2$.
Now, it follows directly from the definition of $A$ and Lemma \ref
{techlmm} that
%
%e33 ###
\begin{eqnarray} \label{eq1}
&& \biggl|\log\sum_{x\in A} e^{hE(x)} - \log\sum_{x\in A} e^{
{\beta
T(x)}/{n} + hE(x)} + \frac{\beta n^2\varphi(u^*)^3}{6}
\biggr|\nonumber\\
&&\qquad= \biggl|\log\sum_{x\in A} e^{hE(x) + {\beta n^2\varphi
(u^*)^3}/{6}} - \log\sum_{x\in A} e^{{\beta T(x)}/{n} + hE(x)}\biggr|
\\
&&\qquad\le\beta\max_{x\in A} \biggl|\frac{T(x)}{n} - \frac{n^2\varphi
(u^*)^3}{6} \biggr|
\le\beta K(\beta,h) n^{3/2}.\nonumber
\end{eqnarray}
Next, observe that
%
%e34 ###
\begin{eqnarray}\label{eq11}
&& \biggl|\log\sum_{x\in A} e^{{\beta T(x)}/{n} + hE(x)} - \log
\sum
_{x\in\Omega} e^{{\beta T(x)}/{n} + hE(x)}
\biggr|\nonumber\\[-8pt]\\[-8pt]
&&\qquad = |{\log}\pr(X\in A) |\le|{\log}(3/4)|.\nonumber
\end{eqnarray}
Similarly, we have
%
%e35 ###
\begin{eqnarray}\label{eq12}
&& \biggl|\log\sum_{x\in B} e^{{\beta T(x)}/{n} + hE(x)} - \log
\sum
_{x\in\Omega} e^{{\beta T(x)}/{n} + hE(x)} \biggr| \nonumber\\[-8pt]\\[-8pt]
&&\qquad= |{\log}\pr(X\in B) |\le|{\log}(1/2)|,\nonumber
\end{eqnarray}
where we used the fact that $\pr(X\in A\cap C)\ge\pr(X\in A)+\pr
(X\in C)-1$.
Combining the last two inequalities, we get
%
%e36 ###
\begin{equation}\label{eq2}
\biggl|\log\sum_{x\in A} e^{{\beta T(x)}/{n} + hE(x)} - \log
\sum
_{x\in B} e^{{\beta T(x)}/{n} + hE(x)} \biggr| \le\log(8/3).
\end{equation}
Next, note that by the definition of $B$ and Lemma \ref{techlmm}, we
have that for any $h'$,
%
%e37 ###
\begin{eqnarray}\label{eq3}
&& \hspace*{20pt}\biggl|\log\sum_{x\in B} e^{{\beta T(x)}/{n} + hE(x)} -\frac
{n^2(h-h')\varphi(u^*)}{2} - \frac{\beta n^2 \varphi(u^*)^3}{6} -
\log
\sum_{x\in B} e^{h'E(x)} \biggr| \nonumber\hspace*{-20pt}\\
&&\hspace*{20pt}\qquad\le\sup_{x\in B} \biggl| \frac{\beta T(x)}{n} + hE(x) -\frac
{n^2(h-h')\varphi(u^*)}{2} - \frac{\beta n^2 \varphi(u^*)^3}{6} -
h'E(x) \biggr|
\\
&&\hspace*{20pt}\qquad\le(\beta+|h-h'|)K(\beta,h)n^{3/2}.\nonumber
\end{eqnarray}
Now, choose $h' = \log\frac{\varphi(u^*)}{1-\varphi(u^*)}$. Then
%
%e38 ###
\begin{equation}\label{eq4}
\biggl|\log\sum_{x\in B} e^{h'E(x)} - \log\sum_{x\in\Omega}
e^{h'E(x)} \biggr| = |{\log\pr}(Z\in B)|\le\log2.
\end{equation}
Adding up (\ref{eq1}), (\ref{eq2}), (\ref{eq3}) and (\ref{eq4}), and
using the triangle inequality, we get
%
%e39 ###
\begin{equation}\label{eq5}\quad
\biggl|\log\sum_{x\in A} e^{hE(x)} - \frac{n^2(h-h')\varphi(u^*)}{2}
-\log\sum_{x\in\Omega}e^{h'E(x)} \biggr|\le K'(\beta, h) n^{3/2},
\end{equation}
where $K'(\beta,h)$ is a constant depending only on $\beta,h$. For any
$s\in\dR$, a trivial verification shows that
\[
\log\sum_{x\in\Omega} e^{sE(x)} = \pmatrix{n\cr2}\log(1+e^s).
\]
Again, note that $\log\pr(Y\in A) = \log\sum_{x\in A} e^{hE(x)} -
\log
\sum_{x\in\Omega} e^{hE(x)}$. Therefore, it follows from inequality
(\ref{eq5}) that
\[
\biggl|\frac{\log\pr(Y\in A)}{n^2} - \frac{(h-h')\varphi(u^*) +
\log
(1+e^{h'}) - \log(1+e^h)}{2} \biggr|\le\frac{K'(\beta, h)}{\sqrt{n}}.
\]
Now $h = \log\frac{\varphi(0)}{1-\varphi(0)}$ and $h' = \log\frac
{\varphi(u^*)}{1-\varphi(u^*)}$. Also, $\log(1+e^h) = -\log
(1-\varphi
(0))$ and $\log(1+e^{h'}) = - \log(1-\varphi(u^*))$. Substituting
these in the above expression, we get
\[
\biggl|\frac{\log\pr(Y\in A)}{n^2} - \frac{-I(\varphi
(u^*),\varphi
(0))}{2} \biggr|\le\frac{K'(\beta,h)}{\sqrt{n}}.
\]
This completes the proof of the lemma.
\end{pf}

We are now ready to finish the proof of Theorem \ref{pressure}.
\begin{pf*}{Proof of Theorem \ref{pressure}}
Note that by adding the terms in (\ref{eq12}), (\ref{eq3}) and
(\ref{eq4}) from the proof of Lemma \ref{large1}, and applying the triangle
inequality, we get
\[
\biggl|\frac{\log Z_n(\beta,h)}{n^2} - \frac{(h-h')\varphi(u)}{2}
-\frac
{\beta\varphi(u)^3}{6} - \frac{1}{2}\log(1+e^{h'}) \biggr|\le\frac
{K(\beta,h)}{\sqrt{n}}.
\]
This can be rewritten as
\[
\biggl|\frac{\log Z_n(\beta,h)}{n^2} + \frac{I(\varphi(u),\varphi
(0)) +
\log(1-\varphi(0))}{2} -\frac{\beta\varphi(u)^3}{6} \biggr|\le
\frac
{K(\beta,h)}{\sqrt{n}}.
\]
%
%The only part yet unproved is the assertion the condition \eqref{cond}
%is satisfied when $h \ge0$ or $\beta< 4$. First, note that if $\beta
%< 4$, then a simple verification shows that
%all } x\in[0,1].
%Condition \eqref{cond} follows directly from this. The case $h\ge0$
%is harder. First, note that if $h \ge0$, then $\varphi(x) \ge
%In other words, $\varphi$ is a concave function on $[0,1]$. Thus, the
%function
%which is always positive because $\varphi$ is increasing, attains its
%maximum in $[0,1]$ at $x=0$. Now $\varphi(u^*) = \sqrt{u^*}$.
%Combining, we get
%(0,1]} \frac{\sqrt{u} - \varphi(0)}{u} = \frac{1}{4\varphi(0)} \le
%
This completes the proof of Theorem \ref{pressure}.
\end{pf*}

Note that the proof of Theorem \ref{pressure} contains a proof for the
lower bound in the general case. We provide the proof below for completeness.

\begin{pf*}{Proof of Lemma \ref{lem:lbd}}
Fix any $r\in(0,1)$. Define the set $B_{r}$ as
\[
B_{r}= \biggl\{x\in\{0,1\}^n\dvtx \biggl|\frac{T(x)}{n} - \frac{n^2r^3}{6}
\biggr|\le K(r) n^{3/2}, \biggl|E(x) - \frac{n^2r}{2} \biggr|\le K(r)
n^{3/2} \biggr\},
\]
where $K(r)$ is chosen in such a way that $\pr(Z\in B_{r})\ge1/2$
where $Z=((Z_{ij}))_{i<j}$ and $Z_{ij}$'s are i.i.d. Bernoulli$(r)$.
From the proof of Lemma \ref{large1}, it is easy to see that
\[
\biggl|\log\sum_{x\in B_r}e^{{\beta T(x)}/{n} + hE(x)} - \frac
{n^{2}}{2} \biggl( (h-h')r +\frac{\beta r^3}{3} +\log(1+e^{h'})
\biggr) \biggr| \le K'n^{3/2},
\]
where $h'=\log\frac{r}{1-r}$ and $K'$ is a constant depending on
$\beta
,h,r$. Simplifying, we have
%
%e40 ###
\begin{eqnarray}\label{eq:lowerbd}
\frac{2}{n^{2}}\log Z_{n}(\beta,h)&\ge&\frac{2}{n^{2}}\log\sum
_{x\in
B_r}e^{{\beta T(x)}/{n} + hE(x)} \nonumber\\[-8pt]\\[-8pt]
&\ge&\frac{\beta r^3}{3} +\log(1-p) -I(r,p) -
\frac{K'}{\sqrt{n}}\nonumber
\end{eqnarray}
for all $r$ where $p=e^{h}/(1+e^{h})$. Now taking limit as $n\to\infty$
and maximizing over $r$ we have the first inequality (\ref{eq:lbd}).
Given $\beta,h$, define the function
\[
f(r)=\frac{\beta r^3}{3} +\log(1-p) -I(r,p),
\]
where $p=e^{h}/(1+e^{h})$. One can easily check that $f'(r)\gtreqqless
0$ iff $\varphi(u)^{2}-u\gtreqqless0$ for $u=r^{2}$. From this fact,
the second equality follows.
\end{pf*}
\begin{lem}\label{large2}
Let $T_n$ be the number of triangles in the Erd\H os--R\'enyi graph
$G(n, \varphi(0))$. Then there is a constant $K(\beta,h)$ depending
only on $\beta$ and $h$ such that for all $n$
\[
\frac{\log\pr(T_n \ge{n \choose3}\varphi(u^*)^3)}{n^2} \le
\frac
{-I(\varphi(u^*),\varphi(0))}{2} + \frac{K(\beta,h)}{\sqrt{n}}.
\]
\end{lem}
\begin{pf}
By Markov's inequality, we have
\[
\frac{\log\pr(T_n \ge{n\choose3}\varphi(u^*)^3)}{n^2} \le
-\frac
{\beta}{n^3}\pmatrix{n \cr3}\varphi(u^*)^3 + \frac{\E(e^{\beta
T_n/n})}{n^2}.
\]
From the last part of Theorem \ref{pressure}, it is easy to obtain an
optimal upper bound of the second term on the right-hand side, which
finishes the proof of the lemma.
\end{pf}
\begin{pf*}{Proof of Theorem \ref{largedev}}
Given $p$ and $r$, if for all $r'$ belonging to a small neighborhood of
$r$ there exist $\beta$ and $h$ satisfying the conditions of Theorem
\ref{pressure} such that $\varphi(0) = p$ and $\varphi(u^*) = r'$, then
a combination of Lemma \ref{large1} and Lemma \ref{large2} implies the
conclusion of Theorem \ref{largedev}. If $p\ge p_{0}=2/(2+e^{3/2})$, we
can just choose $h\ge h_{0}=-\log2 - 3/2$ such that $p = e^h/(1+e^h)$
and conclude, from Theorem~\ref{pressure}, Lemma \ref{large1} and Lemma
\ref{lem:pcond}, that the large deviations limit holds for any $\beta
\ge0$. Varying $\beta$ between $0$ and $\infty$, it is possible to get
for any $r \ge p$ a $\beta$ such that $\varphi(u^*) = r$.

For $p\le p_{0}$, we again choose $h$ such that $\varphi(0) = p$. Note
that $h\le h_{0}$. The large deviations limit should hold for any $r\ge
p$ for which there exists $\beta>0$ such that $r = \varphi(u^*) =
\sqrt
{u^*}$ and $(h,\beta)\in S$. It is not difficult to verify that given
$h$, $u^*$ is a continuously increasing function of $\beta$ in the
regime for which $(h,\beta)\in S$. Recall the settings of Lemma \ref
{lem:pcond}. Thus, the values of $r$ that is allowed is in the set
$(p,p_{*})\cup(p^{*},1]$, where $p^{*},p_{*}$ are the unique
nontouching solutions to the equations
\[
\sqrt{p^{*}}=\frac{e^{\beta_{*}(h)p^{*} +h}}{1+e^{\beta_{*}(h)p^{*}
+h}},\qquad
\sqrt{p_{*}}=\frac{e^{\beta^{*}(h)p_{*} +h}}{1+e^{\beta^{*}(h)p_{*}+h}}.
\]
This completes the proof of Theorem \ref{largedev}.
\end{pf*}

Finally, let us round up by proving Lemma \ref{lem:pcond}.
\begin{pf*}{Proof of Lemma \ref{lem:pcond}}
Fix $h\in\dR$. Define the function
\[
\psi(x;h,\beta):= \varphi(x;h,\beta)^{2}-x,
\]
where
\[
\varphi(x;h,\beta) = \frac{e^{\beta x +h}}{1+e^{\beta x +h}} \qquad\mbox{for }
x\in[0,1].
\]
For simplicity, we will omit $\beta,h$ in $\varphi(x;\beta,h)$ and
$\psi
(x;\beta,h)$ when there is no chance of confusion. Note that $\psi
(0)>0>\psi(1)$. Hence, the equation $\varphi(x;\beta,h)=0$ has at least
one solution. Also we have $\psi'(x) = 2\beta\varphi
(x)^{2}(1-\varphi
(x))-1$ and $\varphi$ is strictly increasing. Hence, the equation
$\psi
'(x) =0$ has at most three solutions. So either the function $\psi$ is
strictly decreasing or there exist two numbers $0<a<b<1$ such that
$\psi
$ is strictly decreasing in $[0,a]\cup[b,1]$ and strictly increasing in
$[a,b]$. From the above observations, it is easy to see that the
equation $\psi(x)=0$ has at most three solutions for any $\beta,h$. If
$\psi(x)=0$ has exactly two solutions, then $\psi'=0$ at one of the solution.

Let $u_{*}=u_{*}(h,\beta)$ and $u^{*}=u^{*}(h,\beta)$ be the smallest and
largest solutions of $\psi(x;h,\beta)=0$, respectively.\vspace*{1pt} If $u_{*}=u^{*}$,
we have a unique solution of $\psi(x)=0$. From the fact that $\frac
{\partial}{\partial\beta}\psi(x;h,\beta)>0$ for all $x\in
[0,1],\beta
\ge0,h\in
\dR$, we can deduce that given $h$, $u_{*}(h,\beta)$ and
$u^{*}(h,\beta)$
are increasing functions of $\beta$. Note that $u_{*}$ is left continuous
and $u^{*}$ is right continuous in $\beta$ given $h$. Also note that
given $h\in\dR$, $u^{*}=u_{*}$ if $\beta>0$ is very small or very large.
So, we can define $\beta_{*}(h)$ and $\beta^{*}(h)$ such that for
$\beta<
\beta
_{*}(h)$ and for $\beta> \beta^{*}(h)$ we have $u_{*}(h,\beta
)=u^{*}(h,\beta)$.
$\beta_{*}$ is the largest and $\beta^{*}$ is the smallest such number.

Therefore, we can deduce that at $\beta=\beta_{*}(h),\beta^{*}(h)$ the
equation $\psi(x;h,\beta)=0$ has exactly two solutions. Thus, we have two
real numbers $x_{*},x^{*}\in[0,1]$ such that
\[
\varphi(x)^{2}=x \quad\mbox{and}\quad 2\beta\varphi(x)^{2}\bigl(1-\varphi(x)\bigr)=1
\]
for $(x,\beta)=(x_{*},\beta_{*})$ or $(x^{*},\beta^{*})$. Thus, we have
$2\beta
x(1-\sqrt{x})=1$ and
\[
h= \log\frac{\sqrt{x}}{1-\sqrt{x}} - \frac{1}{2(1-\sqrt{x})}
\]
for $x=x_{*},x^{*}$. Define $a_{*}= x_{*}^{-1/2}-1$ and $a^{*}=
(x^{*})^{-1/2}-1$. Note that $x=(1+a)^{-2}, \beta=(1+a)^{3}/2a^{2}$ for
$(x,a,\beta)=(x_{*},a_{*},\beta_{*})$ or $(x^{*},a^{*},\beta^{*})$
and we have
%
%e41 ###
\begin{equation}\label{eq:h}
h= -\log a -\frac{1+a}{2a}
\end{equation}
for $a=a_{*},a^{*}$. Now the function $g(x) = -\log x -({1+x})/{2x}$ is
strictly increasing for $x\in(0,1/2]$ and strictly decreasing for
$x\ge
1/2$. So (\ref{eq:h}) has no solution for $h\ge g(1/2)= \log2
-3/2=:h_{0}$. For $h<h_{0}$, equation (\ref{eq:h}) has exactly two
solutions and for $h=h_{0}$ equation (\ref{eq:h}) has one solution. One
can easily check that $\beta_{*}\le\beta^{*}$ implies that $a_{*}\le
a^{*}$. Also from the fact that (\ref{eq:h}) has at most two solutions,
we have that for $\beta\in(\beta_{*},\beta^{*})$ the equation $\psi(u)=0$
has exactly three solutions.
\end{pf*}
%
%%%%%%%%%%%%%%%%%%%%%%%%%%%%%%%%%%%%%

%s3.4 ###
\subsection{\texorpdfstring{Proof of Lemma \protect\ref{lem:boundary}}{Proof of Lemma 13}}

\begin{figure}[b]

\includegraphics{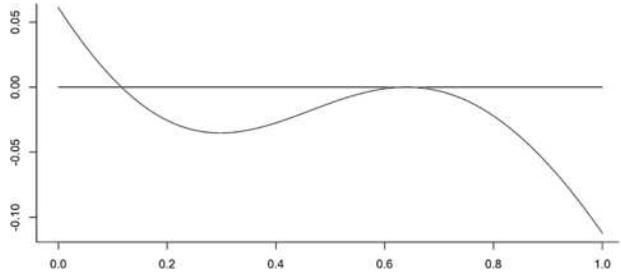}

\caption{The function $\psi(\cdot)$ for $(h,\beta)=\gamma(1/4)$.}
\label{fig:boundary}
\end{figure}

For simplicity, we will prove the result only for the lower boundary
part, that is, for $(h,\beta)=\gamma(t)$ with $t<1/2$. The proof for
the upper boundary is similar. Fix $t<1/2$. Let us briefly recall the
setup. The function $\psi(u)=\varphi(u)^{2}-u$ has two roots at
$0<u^{*}<v^{*}<1$ and $\psi'(u_{*})<0$ while $\psi'(v^{*})=0, \psi
''(v^{*})<0$. See Figure \ref{fig:boundary} for the graph of the
function $\psi$ when $t=1/4$.

%f5 ###

Define the function
\[
f(r)=\frac{\beta r^3}{3} +\log(1-p) -I(r,p) \qquad\mbox{for } r\in(0,1).
\]
From the proof of Lemma \ref{lem:lbd} and the fact that $\psi'(u)< 0$
for $u\in(u^{*},v^{*})$, it is easy to see that $f(\varphi
(u^{*}))>f(\varphi(v^{*}))$ and
%
%e42 ###
\begin{equation}\label{bdlower}
\frac{2}{n^{2}}\log Z_{n}(\beta,h)
\ge f(\varphi(u^{*})) - \frac{K}{\sqrt{n}},
\end{equation}
where $K$ depends on $\beta,h$. Now, using the same idea used in the
proof of Lem\-ma~\ref{mean2}, we have
\[
\pr(
\Delta\ge t ) \le n^{2}\exp\biggl(- \frac{nt^{2}}{8(1+\beta)}
\biggr)
\]
for all $t\ge8\beta/n$ and
$
\psi(L_{\max})\ge-\Delta, \psi(L_{\min})\le\Delta
$
where
\[
\Delta= \max_{1\le i<j\le n} \biggl|L_{ij} - \frac{1}{n}\sum_{k\notin\{
i,j\}} \varphi(L_{ik})\varphi(L_{jk}) \biggr|.
\]
Hence, there exists $\eps_{0}>0,c>0$ such that whenever $\Delta<\eps_{0}$
we have $L_{\min}\ge u^{*}-c\Delta$ and either $L_{\max}\le
u^{*}+c\Delta
$ or
$|L_{\max}-v^{*}|\le c\sqrt{\Delta}$. Define
%
%e43 ###
\begin{equation}
U=\{L_{\max}< (u^{*}+v^{*})/2\}.
\end{equation}
Then again using the idea used in Lemma \ref{mean2} one can easily
show that
\[
\E(\ind_{U}\cdot|L_{ij}-u^{*}|)\le\frac{K(\beta,h)}{n^{1/2}}
\qquad\mbox{for
all }i<j.
\]
We will show that $\pr(U^{c})\le(\log n)^{2}/n$ and it will imply that
\[
\E(|L_{ij}-u^{*}|)\le\E(\ind_{U}\cdot|L_{ij}-u^{*}|) + \pr
(U^{c})\le
\frac{K(\beta,h)}{n^{1/2}} \qquad\mbox{for all }i<j.
\]
Then the rest of the assertions follow using the steps in the proof of
Theorem \ref{large1}.

Hence, let us concentrate on the event $U^{c}$. It is enough to
restrict to the event $U^{c}\cap\{|L_{\max}-v^{*}|\le c\sqrt{\Delta
}\}
\cap
\{L_{\min}\ge u^{*}-c\Delta\}$. Here, the rough idea is that, a~large
fraction of $L_{ij}$'s has to be near $v^{*}$ in order to make $L_{\max
}\simeq v^{*}$.
Suppose $L_{\max}=L_{i_{0}j_{0}}$. Define the set
\[
A=\{k\dvtx L_{i_0k}<L_{\max}-\delta_{1}\},
\]
where $\delta_{1}$ will be chosen later such that $\delta_{1}+c\sqrt
{\Delta
}<v^{*}-u^{*}$. Note that $\varphi(u)^{2}\le\max\{u,u^{*}\}$ for all
$u$ and by assumption $|L_{\max}-v^{*}|\le c\sqrt{\Delta}$. Thus,
$\varphi
(L_{ij})\le\sqrt{L_{\max}}$ for all $i,j$ and $\varphi
(L_{i_{0}k})\le
\sqrt{L_{\max}-\delta_{1}}\le\sqrt{L_{\max}}(1-\delta_{1}/2)$
for $k\in A$.
Thus, we have
\[
L_{\max}=L_{i_{0}j_{0}}\le\Delta+ \frac{1}{n}\sum_{k\neq
i_{0},j_{0}}\varphi(L_{i_{0}k})\varphi(L_{j_{0}k})
\le\Delta+ L_{\max} - \frac{|A|\delta_{1}}{2n},
\]
which clearly implies that
$
\frac{|A|}{n}\le\frac{2\Delta}{\delta_{1}}.
$
Similarly, define the set
$
A_{j}=\{k\dvtx L_{jk}<L_{\max}-\delta_{2}\}
$
where $\delta_{2}$ will be chosen later such that $\delta_{2}+c\sqrt
{\Delta
}<v^{*}-u^{*}$. Using the same idea as before, for $j\notin A$ we have
\[
L_{\max}-\delta_{1}\le L_{i_{0}j}\le\Delta+ L_{\max} - \frac
{|A_{j}|\delta_{2}}{2n}
\quad\mbox{or}\quad \frac{|A_{j}|}{n}\le\frac{2(\Delta+\delta
_{1})}{\delta_{2}}:=M\qquad
(\mbox{say}).
\]
Choose $\delta_{2}=\Delta^{1/5}, \delta_{1}=\Delta^{3/5}$. Then we have
\begin{eqnarray*}
\sum_{i<j}|L_{ij}-L_{\max}|^{2}
&\le&\frac{n|A|+ nM + n^{2}\delta_{2}^{2}}{2} \\
&\le&\frac{n^{2}\Delta}{\delta_{1}} + \frac{n^{2}(\Delta+\delta
_{1})}{\delta_{2}} +
\frac{n^{2}\delta_{2}^{2}}{2} \\
&\le& 4n^{2}\Delta^{2/5}.
\end{eqnarray*}
Thus, by symmetry and H\"{o}lder's inequality, we have
%
%e44 ###
\begin{eqnarray}
\label{vstar}
\E(\ind_{U^{c}}\cdot|L_{ij}-v^{*}|^{2}) &\le& K\E(\ind_{U^{c}}\cdot
\Delta
^{2/5}) \le K\pr(U^{c})^{9/10}\cdot\E(\Delta^{4})^{1/10}
\nonumber\\[-8pt]\\[-8pt]
&\le& \frac{K(\log n)^{1/5}}{n^{1/5}}\pr(U^{c})^{9/10}\nonumber
\end{eqnarray}
for some constant $K$. Now using Lemma \ref{lem:ex} and
(\ref{vstar}) we have
%
%e45 ###
\begin{eqnarray}\label{ucomp}
\E\biggl[\biggl| E(\vX) - \frac{n^{2}\varphi(v^{*})}{2}\biggr|
\Big|U^{c} \biggr]&\le&\frac{Cn^{9/5}(\log
n)^{1/5}}{\pr(U^{c})^{1/10}}\quad\mbox{and}\nonumber\\[-8pt]\\[-8pt]
\E\biggl[\biggl|\frac{T(\vX)}{n} - \frac{n^{2}\varphi
(v^{*})^{3}}{6}\biggr| \Big| U^{c} \biggr] &\le&\frac
{Cn^{9/5}(\log
n)^{1/5}}{\pr(U^{c})^{1/10}}.\nonumber
\end{eqnarray}
If $\pr(U^{c})>(\log n)^{2}/n$, from inequality (\ref{ucomp}) we have
\[
\pr\biggl( \biggl| E(\vX) - \frac{n^{2}\varphi(v^{*})}{2}\biggr|\ge
Kn^{19/10} \big| U^{c} \biggr)\le\frac{1}{4}
\]
and
\[
\pr\biggl( \biggl| \frac{T(\vX)}{n} - \frac
{n^{2}\varphi
(v^{*})^{3}}{6}\biggr|\ge Kn^{19/10} \big| U^{c} \biggr)\le
\frac{1}{4}
\]
for some large constant $K$ depending on $\beta,h$. Now define the set
\[
B= \biggl\{x\in\{0,1\}^n\dvtx \biggl|\frac{T(x)}{n} - \frac{n^2\varphi
(v^{*})^3}{6} \biggr|\le Kn^{19/10}, \biggl|E(x) - \frac{n^2\varphi
(v^{*})}{2} \biggr|\le Kn^{19/10} \biggr\}.
\]
Using the same idea used in the proof of Lemma \ref{large1}, one can
again show that
\[
\biggl|\frac{2}{n^{2}}\log(Z_{n}\pr(U^{c})) - f(\varphi(v^{*}))
\biggr|\le\frac{K}{n^{1/10}}
\]
for some constant $K$ depending on $\beta,h$. The crucial fact is that
$\pr(\{L_{\max}(\vZ)>(u^{*}+v^{*})/2\}\cap\{\vZ\in B\})$ is bounded
away from zero when $\vZ=((Z_{ij}))_{i<j}\sim$ G$(n,\varphi(v^{*}))$.
Thus, we have
\[
\biggl|\frac{2}{n^{2}}\log Z_{n} - f(\varphi(v^{*})) \biggr|\le
\frac
{K}{n^{1/10}}.
\]
But this leads to a contradiction, since by (\ref{bdlower})
we have
\[
\frac{2}{n^{2}}\log Z_{n}(\beta,h)
\ge f(\varphi(u^{*})) - \frac{K}{\sqrt{n}}
\]
and $f(\varphi(u^{*})) > f(\varphi(v^{*})) $. Thus, we have $\pr
(U^{c})\le(\log n)^{2}/n$ and we are done.

%%%%%%%%%%%%%%%%%%%%%%%%%%%%%%%%%%%%%
%s3.5 ###
\subsection{\texorpdfstring{Proof of Theorem
\protect\ref{pressuregen}}{Proof of Theorem 15}}

The proof is almost an exact copy of the proof of Theorem \ref
{pressure}. Recall the definition of $L_{ij}$,
%
%e46 ###
\begin{equation}
L_{ij}:= \frac{N(\vX_{(i,j)}^{1})-N(\vX_{(i,j)}^{0})}{(n-2)_{\mvv
_{F}-2}} \qquad\mbox{for } i<j.
\end{equation}
In fact, we can write $L_{ij}$ explicitly as a horrible sum
\[
L_{ij}= \frac{1}{\alpha_{F}(n-2)_{\mvv_{F}-2}}\mathop{\sum
_{t_{1}<t_{2}<\cdots<t_{\mvv_{F}-2}}}_{t_{l}\in[n]\setminus\{i,j\}
\ \mathrm{for}\ \mathrm{all}\ l} \sum_{(a,b)\in E(F)}{\sum_{\pi}}^{\prime}
\mathop{\prod_{(k,l)\in E(F)}}_{(k,l)\neq(a,b)} X_{\pi_{k}\pi_{l}},
\]
where the sum $\sum'$ is over all one-one onto map $\pi$ from
$V(F)=[\mvv_F]$ to $\{a,b,t_{1},\ldots,t_{\mvv_{F}-2}\}$ where $\{
\pi
(a),\pi(b)\}=\{i,j\}$. Now, we briefly state the main steps. First, we
have $\E(X_{ij}\mid\mbox{rest})=\varphi(L_{ij})$. Moreover, using
Lemma \ref{techlmm} it is easy to see that $|\E( \prod_{j=1}^{k}
X_{i_{2j-1}i_{2j}}\mid\mbox{rest})-\prod_{j=1}^{k}\varphi
(L_{i_{2j-1}i_{2j}})|\le C\beta/n$ for\vspace*{1pt} every distinct pairs
$(i_{1},i_{2}),\ldots,(i_{2k-1},i_{2k})$ where $C$ is an universal constant.

Now, fix $1\le i<j\le n$. Given a configuration $\vX$, construct
another one $\vX'$ in the following way. Choose $\mvv_{F}-2$ distinct
points uniformly at random without replacement from $[n]\setminus\{
i,j\}
$. Replace the coordinates in $\vX$ corresponding to the edges in the
complete subgraph formed by the chosen points including $i,j$ (except
that we do not change $X_{ij}$) by values drawn from the conditional
distribution given the rest of the edges. Call the new configuration
$\vX'$. Define the antisymmetric function
$
F(\vX,\vX'):= (n-2)_{\mvv_{F}-2}(L_{ij}-L'_{ij}).
$
and $f(\vX):=\E(F(\vX,\vX')\mid\vX)$. Using the same idea as before
and Theorem \ref{thm:conc}, we have
%
%e47 ###
\begin{equation}
\pr( \vert L_{ij} - g_{ij} \vert\ge t ) \le\exp\bigl(
-cnt^{2}/(1+\beta)\bigr),
\end{equation}
where $c$ is an absolute constant and $g_{ij}$ is obtained from
$L_{ij}$ by replacing $X_{kl}$ by $\varphi(L_{kl})$ for all $k<l$. Note
that there is a slight difference with the calculation in the triangle
case, since we have to consider collections of edges where some are
modified and some are not. But their contribution will be of the order
of $n^{-1}$. Also the conditions on $\varphi$ arises in the following
way, if all the $L_{ij}$'s are constant, say equal to $u$, then from
the ``mean-field equations'' for $L_{ij}$'s we must have
\begin{eqnarray*}
u&\approx&\frac{1}{\alpha_{F}(n-2)_{\mvv_{F}-2}}\mathop{\sum
_{t_{1}<t_{2}<\cdots<t_{\mvv_{F}-2}}}_{t_{l}\in[n]\setminus\{i,j\}
\ \mathrm{for}\ \mathrm{all}\ l} \sum_{(a,b)\in E(F)}{\sum_{\pi}}' \varphi
(u)^{\mve
_{F}-1}\\
&=& \frac{2\mve_{F}}{\alpha_{F}} \varphi(u)^{\mve_{F}-1}.
\end{eqnarray*}

The next step is to show that under the conditions on $\varphi$, we
have $\E|L_{ij}-u^{*}|\le Kn^{-1/2}$ for all $i<j$ where $K=K(\beta,h)$
is a constant depending only on $\beta,h$. The crucial fact is that the
behavior of the function $\varphi(u)^k-au$ where $a>0$ is a positive
constant and $k\ge2$ is a fixed integer, is same as the behavior of
the function $\varphi(u)^2-u$.

Now it will follow (using the same proof used for Lemma \ref{lem:ex}) that
\[
\E\biggl\vert E(\vX) - \frac{n^{2}\varphi(u^{*})}{2} \biggr\vert\le Cn^{3/2}
\]
and
\[
\E\biggl\vert N(\vX) - \frac{(n)_{\mvv_{F}}\varphi(u^{*})^{\mve
_{F}}}{\alpha_{F}} \biggr\vert\le Cn^{\mvv_{F}-1/2},
\]
where $C$ is a constant depending only on $\beta,h$. The rest of the
proof follows using the arguments used in the proof of Theorem \ref{pressure}.
%%%%%%%%%%%%%%%%%%%%%%%%%%%%%%%%%%%%%%%%%%%%%%%%%%%%%

\begin{pf*}{Proof of Theorem \ref{largedevgen}}
Using the method of proof for the triangle case and the result from
Theorem \ref{pressuregen}, the proof follows easily.
\end{pf*}

%%%%%%%%%%%%%%%%%%%%%%%%%%%%%%%%%%%%%%%%%%%%%%%%%%%%%
%
\begin{pf*}{Proof of Lemma \ref{lem:pcondgen}}
The proof is same as the proof of Lemma \ref{lem:pcond} except for the
constants.
\end{pf*}

%%%%%%%%%%%%%%%%%%%%%%%%%%%%%%%%%%%%%%%%%%%%%%%%%%%%%

%s3.6 ###
\subsection{\texorpdfstring{Proof of Theorem \protect\ref{pr:ising}}{Proof of Theorem 17}}

Suppose $\mvgs$ is drawn from the Gibbs distribution $\mu_{\beta
,h}$. We
construct $\mvgs'$ by taking one step in the heat-bath Glauber dynamics
as follows: choose a position $I$ uniformly at random from $\Omega$, and
replace the $I$th coordinate of $\mvgs$ by an element drawn from the
conditional distribution of the $\sigma_I$ given the rest. It is easy to
see that $(\mvgs,\mvgs')$ is an exchangeable pair. Let
\[
F(\mvgs,\mvgs') := |\Omega|\bigl(m(\mvgs)-m(\mvgs')\bigr)=\sigma_{I}-\sigma'_{I}
\]
be an antisymmetric function in $\mvgs,\mvgs'$. Since the Hamiltonian
is a simple explicit function, one can easily calculate the conditional
distribution of the spin of the particle at position $x$ given the
spins of the rest.
In fact, we have $\E(\sigma_{x}| \{\sigma_{y}, y\neq x\}])= \tanh
(2\beta d
m_{x}(\mvgs))$ where $m_{x}(\mvgs) := \frac{1}{2d}\sum_{y\in
N_{x}}\sigma
_{y}$ is the average spin of the neighbors of $x$ for $x\in\Omega$. Now,
using Fourier--Walsh expansion we can write the function $\tanh(2\beta d
m_{x}(\mvgs))$ as sums of products of spins in the following way. We have
%
%e48 ###
\begin{equation}\label{eq:fwalsh}
\tanh(2d\beta m_{x}(\mvgs)) = \sum_{k=0}^{2d}a_k(\beta)\sum
_{|S|=k,S\subseteq N_x}\sigma_{S},
\end{equation}
where
%
%e49 ###
\begin{equation}
a_{k}(\beta):= \frac{1}{2^{2d}}\sum_{\mvgs\in\{-1,+1\}^{2d}} \tanh
\Biggl(\beta\sum_{i=1}^{2d}\sigma_{i} \Biggr)\prod_{j=1}^{k}\sigma_{j}
\end{equation}
for $k=0,1,\ldots,2d$. It is easy to see that $a_{k}(\beta)=0$ if $k$ is
even and $a_{k}(\beta)$ is a rational function of $\tanh(2\beta)$ if
$k$ is
odd. Note that the dependence of $a_{k}$ on $d$ is not stated
explicitly. Thus, using (\ref{eq:fwalsh}) and the definitions
in (\ref{def:rk}) we have
\begin{eqnarray*}
f(\mvgs) &=& \E[F(\mvgs,\mvgs')|\mvgs] = \frac{1}{|\Omega|}\sum
_{x\in\Omega}
E[ \sigma_{x} -\sigma'_{x} |\mvgs]\\
&=& m(\mvgs) - \frac{1}{|\Omega|}\sum_{x\in\Omega}\tanh(2\beta d
m_{x}(\mvgs))\\
&=& \bigl(1-2da_{1}(\beta)\bigr)m(\mvgs) - \sum_{k=1}^{d-1}\pmatrix{2d\cr
2k+1}a_{2k+1}(\beta) r_{2k+1}(\mvgs).
\end{eqnarray*}
Define $\theta_{k}(\beta):={2d\choose2k+1}a_{2k+1}(\beta)$ for
$k=0,1,\ldots,d-1$.
Note that we can explicitly calculate the value of $\theta_{0}(\beta)$
as follows:
\[
\theta_{0}(\beta)%&= \frac{2d}{4^{d}} \sum_{\mvgs\in\{-1,+1\}^{2d}}
% (\beta\sum_{i=1}^{2d}\sigma_{i} )\sigma_{1}\\
= \frac{1}{4^{d}} \sum_{\mvgs\in\{-1,+1\}^{2d}} \tanh\Biggl(\beta
\sum
_{i=1}^{2d}\sigma_{i} \Biggr)\sum_{i=1}^{2d}\sigma_{i}
= \frac{2}{4^{d}}\sum_{k=1}^{d}2k\pmatrix{2d\cr d+k} \tanh(2k\beta).
\]
Now, we have $|F(\mvgs,\mvgs')|\le2$ and
\[
|f(\mvgs) - f(\mvgs')| \le\frac{2}{|\Omega|} \Biggl( |1-\theta
_{0}(\beta)| +
\sum_{k=1}^{d-1}(2k+1)\theta_{k}(\beta) \Biggr) = \frac{2}{|\Omega|}
b(\beta)
\]
for all values of $\mvgs,\mvgs'$.
Hence, the condition of Theorem \ref{thm:conc} is satisfied with $B=0$,
$C=2|\Omega|^{-1}b(\beta)$.
So by part (ii) of Theorem \ref{thm:conc}, we have
\[
\pr\Biggl( \sqrt{|\Omega|} |\bigl(1-\theta_{0}(\beta)\bigr)m(\mvgs) - \sum
_{k=1}^{d-1}\theta_{k}(\beta) r_{2k+1}(\mvgs) | \ge t
\Biggr)\le
2\exp\biggl( - \frac{t^{2} }{4b(\beta) } \biggr)
\]
for all $t>0$.
Obviously, $\theta_{0}(\cdot)$ is a strictly increasing function of
$\beta
$. Also, we have $\theta_{0}(0)=0$ and
\[
\theta_{0}(\infty):=\lim_{\beta\to\infty}\theta_{0}(\beta) =
\frac
{1}{4^{d-1}}\sum_{k=1}^{d}k\pmatrix{2d\cr d+k}.
\]
For $d=1$, we have $\theta_{0}(\infty)=1$ and for $d\ge2$ we have
\begin{eqnarray*}
\theta_{0}(\infty)&\ge&\frac{1}{4^{d-1}} \Biggl[ 2\sum
_{k=1}^{d}\pmatrix{2d\cr d+k} - \pmatrix{2d\cr d+1} \Biggr]\\
&=& \frac{1}{4^{d-1}} \biggl[ 2^{2d} - \pmatrix{2d\cr d}-\pmatrix{2d\cr
d+1} \biggr]\\
&=& 4 - \frac{8}{2^{2d+1}}\pmatrix{2d+1\cr d+1}
\end{eqnarray*}
and from the fact that $\sum_{k=d-1}^{d+2}{2d+1\choose k}\le2^{2d+1}$
we have
\[
\frac{1}{2^{2d+1}}\pmatrix{2d+1\cr d+1} \le\frac{d+2}{4(d+1)}\le
\frac
{1}{3} \qquad\mbox{for }d\ge2.
\]
Hence, for $d\ge2$ we have $\theta_{0}(\infty)>1$ and there exists
$\beta
_{1}\in(0,\infty)$, depending on $d$, such that $1-\theta_{0}(\beta)>0$
for $\beta<\beta_{1}$ and $1-\theta_{0}(\beta)<0$ for $\beta>\beta
_{1}$. This
completes the proof.
%%%%%%%%%%%%%%%%%%%%%%%%%%%%%%%%%%%%%%%

%s3.7 ###
\subsection{\texorpdfstring{Proof of Proposition
\protect\ref{pr:isingh}}{Proof of Proposition 19}}

The proof is almost same as the proof of Proposition \ref{pr:ising}.
Define $\mvgs,\mvgs'$ as before. Define the antisymmetric function
$F(\mvgs,\mvgs')$ as follows:
\begin{eqnarray*}
F(\mvgs,\mvgs') :\!&=& |\Omega|\bigl(1+\tanh(h)\tanh(2\beta d m_{I}(\mvgs
))\bigr)\bigl(m(\mvgs
)-m(\mvgs')\bigr)\\
&=&\bigl(1+\tanh(h)\tanh(2\beta d m_{I}(\mvgs))\bigr)(\sigma_{I}-\sigma'_{I}).
\end{eqnarray*}
Recall that $m_{x}(\mvgs) := \frac{1}{2d}\sum_{y\in N_{x}}\sigma
_{y}$ is
the average spin of the neighbors of $x$ for $x\in\Omega$. Now under
$\mu
_{\beta,h}$, we have
\begin{eqnarray*}
\E(\sigma_{x}| \{\sigma_{y}, y\neq x\}) &=& \tanh\bigl(2\beta d
m_{x}(\mvgs)+h\bigr)\\
&=&\frac{\tanh(h)+ \tanh(2\beta d m_{x}(\mvgs))}{1+\tanh(h)\tanh
(2\beta d
m_{x}(\mvgs))}.
\end{eqnarray*}
Thus, we have
\begin{eqnarray*}
f(\mvgs) &=& \E(F(\mvgs,\mvgs')|\mvgs) \\
&=& \frac{1}{|\Omega|}\sum_{x\in\Omega}\bigl(1+\tanh(h)\tanh(2\beta d
m_{x}(\mvgs))\bigr)
\E( \sigma_{x} -\sigma'_{x} |\mvgs)\\
&=& m(\mvgs) -\tanh(h)+ \frac{1}{|\Omega|}\sum_{x\in\Omega}\bigl(\tanh
(h)\sigma
_{x}-1\bigr)\tanh(2\beta d m_{x}(\mvgs)).
\end{eqnarray*}
After some simplifications and using the definitions of the functions,
$r,s$ we have
\begin{eqnarray*}
f(\mvgs)&=& \bigl(1-\theta_{0}(\beta)\bigr)m(\mvgs) - \sum_{k=1}^{d-1}\theta
_{k}(\beta
) r_{2k+1}(\mvgs)\\
&&{} -\tanh(h) \Biggl(1-\sum_{k=0}^{d-1}\theta_{k}(\beta)
s_{2k+1}(\mvgs) \Biggr).
\end{eqnarray*}
Now for all values of $\mvgs,\mvgs'$ we have
\[
|f(\mvgs) - f(\mvgs')|\le\frac{2}{|\Omega|} b(\beta)(1+\tanh|h|)
\]
and the proof henceforth is exactly as in the proof of Proposition \ref{pr:ising}.
%%%%%%%%%%%%%%%%%%%%%%%%%%%%%%%%%%%%%%%%%%%%%%%%%

%s3.8 ###
\subsection{\texorpdfstring{Proof of Theorem \protect\ref{thm:nong1}}{Proof of Theorem 2}}

Assume that $\psi(0)>0$. We will handle the case $\psi(0)=0$ later.
Note that condition (\ref{eq:al}) implies that $x^\alpha/\psi(x)$ is a
nondecreasing function for $x>0$. Define the function
\[
\varphi(x):=\frac{x^2}{\psi(x)} \quad\mbox{and}\quad \gamma(x):= 2-\frac
{x\psi
'(x)}{\psi(x)} \qquad\mbox{for } x\neq0
\]
and $\varphi(0)=0,\gamma(0)=2$.
Clearly, we have $2-\alpha\le\gamma(x)\le2$ for all $x\in\dR$.
Now, $\limsup_{x\to0}\varphi(x)\le\lim_{x\to0+}x^{2-\alpha}/\psi
(1)=0=\varphi(0)$ as $\alpha<2$. Also $\varphi(x)$ is differentiable in
$\dR\setminus\{0\}$ with
%
%e50 ###
\begin{equation}\label{eq:lphi}
\varphi'(x)= \frac{x \gamma(x)}{\psi(x)} > 0 \qquad\mbox{for } x\neq0.
\end{equation}
Hence, $\varphi$ is absolutely continuous in $\dR$ and is increasing
for $x\ge0$.

Define $Y=f(X)$. First, we will prove that all moments of $\varphi(Y)$
are finite. Next, we will estimate the moments which will in turn show
that $\varphi(Y)^{1/2}$ has finite exponential moment in $\dR$.
Finally, using Chebyshev's inequality we will prove the tail probability.

By monotonicity of $\psi$ in $[0,\infty)$ and definition of $\alpha$,
we have
%
%e51 ###
\begin{equation}
0\le\frac{x\psi'(x)}{\psi(x)} \le\alpha\qquad\mbox{for all } x\ge0.
\end{equation}
It also follows from (\ref{eq:lphi}) that $0\le(\log\varphi(x))'\le
2/x$ for $x>0$ and integrating we have $\varphi(x)\le\varphi(1)x^2$
for all $x\ge1$. Hence, $\varphi(x)=\varphi(|x|)\le\varphi(1)(1+x^2)$
for all $x\in\dR$ and this, combined with our assumption that $\E
(|f(X)|^{k}) < \infty$ for all $k\ge1$, implies that
$
\E(\varphi(Y)^k)<\infty\mbox{ for all } k\ge1.
$

Define
\[
\beta:= \biggl\lceil\frac{5(2-\alpha)+\delta+1/4}{(2-\alpha)^2}
\biggr\rceil\ge3.
\]
Fix an integer $k\ge\beta$ and define
\[
g(x)=\frac{x^{2k-1}}{\psi^k(x)} \quad\mbox{and}\quad h(x)=\frac
{x^{2k-2}}{\psi
^k(x)} \qquad\mbox{for } x\in\dR.
\]
Clearly, $\E(|Yg(Y)|)<\infty$. Note that $g,h$ are continuously
differentiable in $\dR$ as $k\ge3$. Moreover, for $x\in\dR$ we have,
$|g'(x)|= h(x) \vert k\gamma(x)-1 \vert\le(2k-1)h(x)$, $ h'(x) =
(k\gamma(x)
-2 ) {x^{2k-3}}/{\psi^k(x)}$ and
\[
h''(x)= \bigl[ \bigl(k\gamma(x) -2 \bigr) \bigl(k\gamma(x)
-3 \bigr) + kx\gamma
'(x)\bigr]\frac{x^{2k-4}}{\psi^k(x)}.
\]
We also have
\[
x\gamma'(x)=-\frac{x\psi'(x)}{\psi(x)} \biggl(1-\frac{x\psi
'(x)}{\psi
(x)} \biggr)-\frac{x\psi''(x)}{\psi(x)}\ge-1/4-\delta
\]
for $x\in\dR$. Now $k\ge\beta$ implies that
\begin{eqnarray*}
&& \bigl(k\gamma(x) -2 \bigr) \bigl(k\gamma(x) -3 \bigr) + kx\gamma
'(x)\\
&&\qquad
\ge\bigl(k(2-\alpha) -2 \bigr) \bigl(k(2-\alpha) -3 \bigr) - k(\delta
+1/4)\ge0
\end{eqnarray*}
for all $x$.
Thus, $h''(x)\ge0$ for all $x$ and $h$ is convex in $\dR$.

%Now,
% h'(x) &= (k\gc(x) -2 ) \frac{x^{2k-3}}{\psi^k(x)}
%for $x\ge0$ and $h'(-x)=-h'(-x)$ for $x\le0$.
%Thus $h$ is nonincreasing for $x\le0$ and nondecreasing for $x\ge0$.
%From this fact it is easy to see that
%for any $a\le b$.

Let $X',F(X,X')$ be as given in the hypothesis. Define $Y'=f(X')$.
Recall that $(X,X')$ is an exchangeable pair and so is $(Y,Y')$. Using
the fact that $f(X)=\E(F(X,X')|X)$ almost surely, exchangeability of
$(X,X')$ and antisymmetry of $F$, we have
%
%e52 ###
\begin{eqnarray} \label{eq:ln1}
\E(Yg(Y))&=&\E(f(X)g(Y))=\E(F(X,X')g(Y))\nonumber\\[-8pt]\\[-8pt]
&=&\tfrac{1}{2} \E\bigl(F(X,X')\bigl(g(Y)-g(Y')\bigr) \bigr).\nonumber
\end{eqnarray}
Now, for any $x< y$ we have
\[
\biggl\vert\frac{g(x)-g(y)}{x-y} \biggr\vert=\biggl\vert\int_0^1 g'\bigl(tx+(1-t)y\bigr) \,dt
\biggr\vert\le(2k-1)
\int_0^1 h\bigl(tx+(1-t)y\bigr) \,dt
\]
and convexity of $h$ implies that
\[
\int_0^1 h\bigl(tx+(1-t)y\bigr) \,dt \le\int_0^1 \bigl(t h(x)+(1-t)h(y)\bigr) \,dt =\bigl(h(x)+h(y)\bigr)/2.
\]
Hence, from (\ref{eq:ln1}), we have
%
%e53 ###
\begin{eqnarray} \label{eq:gdh}
\E(Yg(Y))&\le&\frac{2k-1}{4} \E\bigl(|(Y-Y')F(X,X')|\bigl(h(Y)+h(Y')\bigr)\bigr)
\nonumber\\[-8pt]\\[-8pt]
&=&(2k-1)\E(\Delta(X)h(Y)) \le(2k-1)\E(\psi(Y)h(Y)),\nonumber
\end{eqnarray}
where the equality follows by definition of $\Delta(X)$ and
exchangeability of $(Y,Y')$. Thus, for any $k\ge\beta$ we have, from
(\ref{eq:gdh}),
%
%e54 ###
\begin{equation} \label{eq:rec}
\E(\varphi(Y)^k) \le
(2k-1)\E(\varphi(Y)^{k-1}).
\end{equation}
Using induction for $k\ge\beta$, we have
\[
\E(\varphi(Y)^k) \le\frac{(2k)! 2^\beta\beta!}{2^k k!(2\beta
)!}\E
(\varphi
(Y)^{\beta}) \qquad\mbox{for } k\ge\beta.
\]
Also H\"older's inequality applied to (\ref{eq:rec}) for $k=\beta$
implies that $\E(\varphi(Y)^{\beta}) \le(2\beta-1)^\beta$. Thus,
we have
%
%e55 ###
\begin{equation}\label{eq:rec2}
\E(\varphi(Y)^k)\le\cases{
\dfrac{(2k)! 2^\beta\beta!}{k! 2^k(2\beta)!}\E(\varphi(Y)^{\beta
}), &\quad if $k>
\beta$,\vspace*{2pt}\cr
(2\beta-1)^k, &\quad if $0\le k\le\beta$.}
\end{equation}
Note that we have $e^{x}\le e^{x}+e^{-x}=2\sum_{k\ge0}x^{2k}/(2k)!$
for all $x\in\dR$. Combining everything, we finally have
\begin{eqnarray*}
\E(\exp(\theta\varphi(Y)^{1/2}))
&\le&2\sum_{k=0}^\infty\frac{\theta^{2k} }{(2k)!}\E(\varphi
(Y)^{k})\\
&\le&\frac{2^{\beta+1} \beta!}{(2\beta)!}\E(\varphi(Y)^{\beta})
\sum
_{k=\beta
}^\infty\frac{\theta^{2k} }{2^k k!} +
\sum_{k=0}^{\beta-1} \frac{2(2\beta-1)^k\theta^{2k}}{(2k)!} \\
&\le& C_{\beta}\exp(\theta^{2}/2)
\end{eqnarray*}
for all $\theta\ge0$ where the constant $C_{\beta}$ is given by
\[
C_{\beta}:= \max\biggl\{ \frac{ 2(2\beta-1)^{k} 2^{k} k! }{(2k)!} \Big|
0\le k \le\beta\biggr\}.
\]
Here, we used the fact that $(2k)!\ge2^{2k-1}k!^2/k$. Now recall that
$\varphi$ is an increasing function in $[0,\infty)$. Thus, using
Chebyshev's inequality for $\exp(\theta\varphi(x)^{1/2})$ with
$\theta
=\varphi(t)^{1/2}$ we have
\[
\pr\bigl(|f(X)|\ge t\bigr)\le C_{\beta}e^{-\theta\varphi(t)^{1/2}+\theta^{2}/2}=
C_{\beta}e^{-\varphi(t)/2}.
\]

Now suppose that $\psi(0)=0$. For $\eps>0$ fixed, define $\psi_\eps
(x)=\psi(x)+\eps$. Clearly, we have $\Delta(X)\le\psi_{\eps}(f(X))$
a.s. and $\psi_{\eps}$ satisfies all the other properties of $\psi$ including
\[
x\psi'_{\eps}(x)/\psi_{\eps}(x) = x\psi'(x)/\psi(x) \cdot{\psi
(x)}/{\bigl(\psi(x)+\eps\bigr)}\le\alpha
\]
and
\[
x\psi''_{\eps}(x)/\psi_{\eps}(x) = x\psi''(x)/\psi(x)
\cdot{\psi(x)}/{\bigl(\psi(x)+\eps\bigr)}\le\delta
\]
for all $x>0$.
Hence, all the above results hold for $\psi_\eps$ and $\varphi_\eps
(x)=x^{2}/\psi_{\eps}(x)$. Now $\varphi_\eps\uparrow\varphi$ as
$\eps
\downarrow0$. Letting $\eps\downarrow0$, we have the result.

When $\psi$ is once differentiable with $\alpha<2$, it is easy to see that
the function $h$ is nondecreasing (need not be convex) in $[0,\infty)$
for $k\ge\beta:= \lceil2/(2-\alpha)\rceil$. In that case, we have
\[
\int_{0}^{1}h\bigl(tx+(1-t)y\bigr)\,dy\le\max_{z\in[x,y]}h(z)\le h(x)+h(y)
\]
for $x\le y$. Hence, we have the recursion
%
%e56 ###
\begin{equation}
\E(\varphi(Y)^k) \le
2(2k-1)\E(\varphi(Y)^{k-1})
\end{equation}
for $k\ge\beta$. Using the same proof as before, it then follows that
\[
\pr\bigl(|f(X)|\ge t\bigr)\le Ce^{-\varphi(t)/4},
\]
where $C$ depends only on $\alpha$.

%%%%%%%%%%%%%%%%%%%%%%%%%%%%%%%%%%%%%%%%%%%%%%%

\section*{Acknowledgments}
The authors thank Amir Dembo, Erwin Bolthausen, Ofer Zeitouni and Persi
Diaconis for various helpful discussions and comments. They would also
like to thank an anonymous referee for a careful reading of this
article and for constructive criticism that resulted in an improved
exposition.
% helpful comments, pointing out several typos and encouraging them to
%write down Theorem \ref{thm:nong1} in the most general form.

% imsref loaded by lrinkeviciute, 2010-07-01 10:57:58
%

%
\printaddresses


\begin{thebibliography}{39}

%b1 ###
\bibitem{bcr05}
%
\begin{barticle}[mr]
\bauthor{\bsnm{Barthe},~\bfnm{F.}\binits{F.}},
\bauthor{\bsnm{Cattiaux},~\bfnm{P.}\binits{P.}} \AND
\bauthor{\bsnm{Roberto},~\bfnm{C.}\binits{C.}}
(\byear{2005}).
\btitle{Concentration for independent random variables with heavy tails}.
\bjournal{AMRX Appl. Math. Res. Express}
\bvolume{2}
\bpages{39--60}.
\bid{mr={2173316}}
\end{barticle}
%
\endbibitem

%b2 ###
\bibitem{bcr06}
%
\begin{barticle}[vtex]
\bauthor{\bsnm{Barthe},~\bfnm{Franck}\binits{F.}},
\bauthor{\bsnm{Cattiaux},~\bfnm{Patrick}\binits{P.}} \AND
\bauthor{\bsnm{Roberto},~\bfnm{Cyril}\binits{C.}}
(\byear{2006}).
\btitle{Interpolated inequalities between exponential and {G}aussian, {O}rlicz
hypercontractivity and isoperimetry}.
\bjournal{Rev. Mat. Iberoamericana}
\bvolume{22}
\bpages{993--1067}.
\bid{mr={2320410}}
\end{barticle}
%
\endbibitem

%b3 ###
\bibitem{bhamidi08}
%
\begin{bincollection}[vtex]
\bauthor{\bsnm{Bhamidi},~\bfnm{S.}\binits{S.}},
\bauthor{\bsnm{Bresler},~\bfnm{G.}\binits{G.}} \AND
\bauthor{\bsnm{Sly},~\bfnm{A.}\binits{A.}}
(\byear{2008}).
\btitle{Mixing time of exponential random graphs}. In
\bbooktitle{Proc. of the 49th Annual IEEE Symp. on
FOCS}
\bpages{803--812}.
\bpublisher{IEEE Computer Society},
\baddress{Washington, DC}.
\end{bincollection}
%
\endbibitem

%b4 ###
\bibitem{bob07}
%
\begin{barticle}[mr]
\bauthor{\bsnm{Bobkov},~\bfnm{Sergey~G.}\binits{S.~G.}}
(\byear{2007}).
\btitle{Large deviations and isoperimetry over convex probability
measures with
heavy tails}.
\bjournal{Electron. J. Probab.}
\bvolume{12}
\bpages{1072--1100 (electronic)}.
\bid{mr={2336600}}
\end{barticle}
%
\endbibitem

%b5 ###
\bibitem{bl00}
%
\begin{barticle}[vtex]
\bauthor{\bsnm{Bobkov},~\bfnm{S.~G.}\binits{S.~G.}} \AND
\bauthor{\bsnm{Ledoux},~\bfnm{M.}\binits{M.}}
(\byear{2000}).
\btitle{From {B}runn--{M}inkowski to {B}rascamp--{L}ieb and to logarithmic
{S}obolev inequalities}.
\bjournal{Geom. Funct. Anal.}
\bvolume{10}
\bpages{1028--1052}.
\bid{doi={10.1007/PL00001645}, mr={1800062}}
\end{barticle}
%
\endbibitem

%b6 ###
\bibitem{bollobas85}
%
\begin{bbook}[mr]
\bauthor{\bsnm{Bollob{\'a}s},~\bfnm{B{\'e}la}\binits{B.}}
(\byear{2001}).
\btitle{Random Graphs},
\bedition{2nd} ed.
\bseries{Cambridge Studies in Advanced Mathematics}
\bvolume{73}.
\bpublisher{Cambridge Univ. Press}, \baddress{Cambridge}.
\bid{mr={1864966}}
\end{bbook}
%
\endbibitem

%b7 ###
\bibitem{bolthausen87}
%
\begin{barticle}[mr]
\bauthor{\bsnm{Bolthausen},~\bfnm{E.}\binits{E.}}
(\byear{1987}).
\btitle{Laplace approximations for sums of independent random vectors. {II}.
{D}egenerate maxima and manifolds of maxima}.
\bjournal{Probab. Theory Related Fields}
\bvolume{76}
\bpages{167--206}.
\bid{doi={10.1007/BF00319983}, mr={906774}}
\end{barticle}
%
\endbibitem

%b8 ###
\bibitem{bolthausenetal09}
%
\begin{bmisc}[vtex]
\bauthor{\bsnm{Bolthausen},~\bfnm{E.}\binits{E.}},
\bauthor{\bsnm{Comets},~\bfnm{F.}\binits{F.}} \AND
\bauthor{\bsnm{Dembo},~\bfnm{A.}\binits{A.}}
(\byear{2009}).
\bhowpublished{Large
deviations for random matrices and random graphs. Unpublished manuscript}.
\end{bmisc}
%
\endbibitem

%b9 ###
\bibitem{blm03}
%
\begin{barticle}[mr]
\bauthor{\bsnm{Boucheron},~\bfnm{St{\'e}phane}\binits{S.}},
\bauthor{\bsnm{Lugosi},~\bfnm{G{\'a}bor}\binits{G.}} \AND
\bauthor{\bsnm{Massart},~\bfnm{Pascal}\binits{P.}}
(\byear{2003}).
\btitle{Concentration inequalities using the entropy method}.
\bjournal{Ann. Probab.}
\bvolume{31}
\bpages{1583--1614}.
\bid{doi={10.1214/aop/1055425791}, mr={1989444}}
\end{barticle}
%
\endbibitem

%b10 ###
\bibitem{chatterjee05}
%
\begin{bmisc}[vtex]
\bauthor{\bsnm{Chatterjee},~\bfnm{Sourav}\binits{S.}}
(\byear{2005}).
\btitle{Concentration inequalities with exchangeable pairs}.
\bnote{Ph.D. thesis, Stanford Univ. Available at}
\href{http://arxiv.org/abs/arXiv:math/0507526}{arXiv:math/0507526}.
\bid{mr={2622431}}
\end{bmisc}
%
\endbibitem

%b11 ###
\bibitem{chatterjee07}
%
\begin{barticle}[mr]
\bauthor{\bsnm{Chatterjee},~\bfnm{Sourav}\binits{S.}}
(\byear{2007}).
\btitle{Stein's method for concentration inequalities}.
\bjournal{Probab. Theory Related Fields}
\bvolume{138}
\bpages{305--321}.
\bid{doi={10.1007/s00440-006-0029-y}, mr={2288072}}
\end{barticle}
%
\endbibitem

%b12 ###
\bibitem{chatterjee07a}
%
\begin{barticle}[mr]
\bauthor{\bsnm{Chatterjee},~\bfnm{Sourav}\binits{S.}}
(\byear{2007}).
\btitle{Concentration of {H}aar measures, with an application to random
matrices}.
\bjournal{J. Funct. Anal.}
\bvolume{245}
\bpages{379--389}.
\bid{doi={10.1016/j.jfa.2007.01.003}, mr={2309833}}
\end{barticle}
%
\endbibitem

%b13 ###
\bibitem{chatterjeeshao09}
%
\begin{bmisc}[vtex]
\bauthor{\bsnm{Chatterjee},~\bfnm{S.}\binits{S.}} \AND
\bauthor{\bsnm{Shao},~\bfnm{Q.-M.}\binits{Q.-M.}}
(\byear{2009}).
\bhowpublished{Stein's method of exchangeable
pairs with application to the Curie--Weiss model. Preprint. Available
at}
\href{http://arxiv.org/abs/arXiv:0907.4450}{arXiv:0907.4450}.
\end{bmisc}
%
\endbibitem

%b14 ###
\bibitem{chazottes06}
%
\begin{barticle}[mr]
\bauthor{\bsnm{Chazottes},~\bfnm{J.~R.}\binits{J.~R.}},
\bauthor{\bsnm{Collet},~\bfnm{P.}\binits{P.}},
\bauthor{\bsnm{K{\"u}lske},~\bfnm{C.}\binits{C.}} \AND
\bauthor{\bsnm{Redig},~\bfnm{F.}\binits{F.}}
(\byear{2007}).
\btitle{Concentration inequalities for random fields via coupling}.
\bjournal{Probab. Theory Related Fields}
\bvolume{137}
\bpages{201--225}.
\bid{doi={10.1007/s00440-006-0026-1}, mr={2278456}}
\end{barticle}
%
\endbibitem

%b15 ###
\bibitem{doringeichelsbacher09}
%
\begin{barticle}[mr]
\bauthor{\bsnm{D{\"o}ring},~\bfnm{Hanna}\binits{H.}} \AND
\bauthor{\bsnm{Eichelsbacher},~\bfnm{Peter}\binits{P.}}
(\byear{2009}).
\btitle{Moderate deviations in a random graph and for the spectrum of
{B}ernoulli random matrices}.
\bjournal{Electron. J. Probab.}
\bvolume{14}
\bpages{2636--2656}.
\bid{mr={2570014}}
\end{barticle}
%
\endbibitem

%b16 ###
\bibitem{eichelsbacherlowe09}
%
\begin{bmisc}[vtex]
\bauthor{\bsnm{Eichelsbacher},~\bfnm{P.}\binits{P.}} \AND
\bauthor{\bsnm{Lowe},~\bfnm{M.}\binits{M.}}
(\byear{2009}).
\bhowpublished{Stein's method for
dependent random variables occurring in statistical mechanics. Preprint. Available
at}
\href{http://arxiv.org/abs/arXiv:0908.1909}{arXiv:0908.1909}.
\end{bmisc}
%
\endbibitem

%b17 ###
\bibitem{ellis78}
%
\begin{barticle}[vtex]
\bauthor{\bsnm{Ellis},~\bfnm{Richard~S.}\binits{R.~S.}} \AND
\bauthor{\bsnm{Newman},~\bfnm{Charles~M.}\binits{C.~M.}}
(\byear{1978}).
\btitle{The statistics of {C}urie--{W}eiss models}.
\bjournal{J. Stat. Phys.}
\bvolume{19}
\bpages{149--161}.
\bid{mr={0503332}}
\end{barticle}
%
\endbibitem

%b18 ###
\bibitem{ellis78a}
%
\begin{barticle}[mr]
\bauthor{\bsnm{Ellis},~\bfnm{Richard~S.}\binits{R.~S.}} \AND
\bauthor{\bsnm{Newman},~\bfnm{Charles~M.}\binits{C.~M.}}
(\byear{1978}).
\btitle{Limit theorems for sums of dependent random variables
occurring in
statistical mechanics}.
\bjournal{Z. Wahrsch. Verw. Gebiete}
\bvolume{44}
\bpages{117--139}.
\bid{mr={0503333}}
\end{barticle}
%
\endbibitem

%b19 ###
\bibitem{ellis85}
%
\begin{bbook}[mr]
\bauthor{\bsnm{Ellis},~\bfnm{Richard~S.}\binits{R.~S.}}
(\byear{1985}).
\btitle{Entropy, Large Deviations, and Statistical Mechanics}.
\bseries{Grundlehren der Mathematischen Wissenschaften [Fundamental Principles
of Mathematical Sciences]}
\bvolume{271}.
\bpublisher{Springer}, \baddress{New York}.
\bid{mr={793553}}
\end{bbook}
%
\endbibitem

%b20 ###
\bibitem{ggm05}
%
\begin{barticle}[mr]
\bauthor{\bsnm{Gentil},~\bfnm{Ivan}\binits{I.}},
\bauthor{\bsnm{Guillin},~\bfnm{Arnaud}\binits{A.}} \AND
\bauthor{\bsnm{Miclo},~\bfnm{Laurent}\binits{L.}}
(\byear{2005}).
\btitle{Modified logarithmic {S}obolev inequalities and transportation
inequalities}.
\bjournal{Probab. Theory Related Fields}
\bvolume{133}
\bpages{409--436}.
\bid{doi={10.1007/s00440-005-0432-9}, mr={2198019}}%
\end{barticle}%
%
\endbibitem%

%b21 ###
\bibitem{gozlan07}
%
\begin{barticle}[mr]
\bauthor{\bsnm{Gozlan},~\bfnm{Nathael}\binits{N.}}
(\byear{2007}).
\btitle{Characterization of {T}alagrand's like transportation-cost inequalities
on the real line}.
\bjournal{J. Funct. Anal.}
\bvolume{250}
\bpages{400--425}.
\bid{doi={10.1016/j.jfa.2007.05.025}, mr={2352486}}
\end{barticle}
%
\endbibitem

%b22 ###
\bibitem{gozlan09}
%
\begin{barticle}[mr]
\bauthor{\bsnm{Gozlan},~\bfnm{Nathael}\binits{N.}}
(\byear{2010}).
\btitle{Poincare inequalities and dimension free concentration of measure}.
\bjournal{Ann. Inst. H. Poincar\'e Probab. Statist.}
\bnote{To appear.}
\end{barticle}
%
\endbibitem

%b24 ###
\bibitem{ising25}
%
\begin{barticle}[vtex]
\bauthor{\bsnm{Ising},~\bfnm{E.}\binits{E.}}
(\byear{1925}).
\btitle{Beitrag zur theorie des ferromagnetismus}.
\bjournal{Zeitschrift f{\"u}r Physik A Hadrons and Nuclei}
\bvolume{31}
\bpages{253--258}.
\end{barticle}
%
\endbibitem

%b25 ###
\bibitem{jansonetal00}
%
\begin{bbook}[mr]
\bauthor{\bsnm{Janson},~\bfnm{Svante}\binits{S.}},
\bauthor{\bsnm{{\L}uczak},~\bfnm{Tomasz}\binits{T.}} \AND
\bauthor{\bsnm{Rucinski},~\bfnm{Andrzej}\binits{A.}}
(\byear{2000}).
\btitle{Random Graphs}.
\bpublisher{Wiley}, \baddress{New York}.
\bid{mr={1782847}}
\end{bbook}
%
\endbibitem

%b26 ###
\bibitem{jansonolesz04}
%
\begin{barticle}[mr]
\bauthor{\bsnm{Janson},~\bfnm{Svante}\binits{S.}},
\bauthor{\bsnm{Oleszkiewicz},~\bfnm{Krzysztof}\binits{K.}} \AND
\bauthor{\bsnm{Ruci{\'n}ski},~\bfnm{Andrzej}\binits{A.}}
(\byear{2004}).
\btitle{Upper tails for subgraph counts in random graphs}.
\bjournal{Israel J. Math.}
\bvolume{142}
\bpages{61--92}.
\bid{doi={10.1007/BF02771528}, mr={2085711}}
\end{barticle}
%
\endbibitem

%b27 ###
\bibitem{jansonrucinski02}
%
\begin{barticle}[vtex]
\bauthor{\bsnm{Janson},~\bfnm{Svante}\binits{S.}} \AND
\bauthor{\bsnm{Ruci{\'n}ski},~\bfnm{Andrzej}\binits{A.}}
(\byear{2002}).
\btitle{The infamous upper tail: Probabilistic methods in combinatorial optimization}.
\bjournal{Random Structures Algorithms}
\bvolume{20}
\bpages{317--342}.
\bid{doi={10.1002/rsa.10031}, mr={1900611}}
\end{barticle}
%
\endbibitem

%b28 ###
\bibitem{kimvu04}
%
\begin{barticle}[mr]
\bauthor{\bsnm{Kim},~\bfnm{J.~H.}\binits{J.~H.}} \AND
\bauthor{\bsnm{Vu},~\bfnm{V.~H.}\binits{V.~H.}}
(\byear{2004}).
\btitle{Divide and conquer martingales and the number of triangles in
a random
graph}.
\bjournal{Random Structures Algorithms}
\bvolume{24}
\bpages{166--174}.
\bid{doi={10.1002/rsa.10113}, mr={2035874}}
\end{barticle}
%
\endbibitem

%b29 ###
\bibitem{lo00}
%
\begin{bincollection}[mr]
\bauthor{\bsnm{Lata{\l}a},~\bfnm{R.}\binits{R.}} \AND
\bauthor{\bsnm{Oleszkiewicz},~\bfnm{K.}\binits{K.}}
(\byear{2000}).
\btitle{Between {S}obolev and {P}oincar\'e}.
In \bbooktitle{Geometric Aspects of Functional Analysis}.
\bseries{Lecture Notes in Math.}
\bvolume{1745}
\bpages{147--168}.
\bpublisher{Springer}, \baddress{Berlin}.
\bid{doi={10.1007/BFb0107213}, mr={1796718}}
\end{bincollection}
%
\endbibitem

%b30 ###
\bibitem{ledoux01}
%
\begin{bbook}[mr]
\bauthor{\bsnm{Ledoux},~\bfnm{Michel}\binits{M.}}
(\byear{2001}).
\btitle{The Concentration of Measure Phenomenon}.
\bseries{Mathematical Surveys and Monographs}
\bvolume{89}.
\bpublisher{Amer. Math. Soc.}, \baddress{Providence, RI}.
\bid{mr={1849347}}
\end{bbook}
%
\endbibitem

%b31 ###
\bibitem{martinlof82}
%
\begin{barticle}[mr]
\bauthor{\bsnm{Martin-L{\"o}f},~\bfnm{Anders}\binits{A.}}
(\byear{1982}).
\btitle{A {L}aplace approximation for sums of independent random variables}.
\bjournal{Z. Wahrsch. Verw. Gebiete}
\bvolume{59}
\bpages{101--115}.
\bid{doi={10.1007/BF00575528}, mr={643791}}
\end{barticle}
%
\endbibitem

%b32 ###
\bibitem{onsager44}
%
\begin{barticle}[mr]
\bauthor{\bsnm{Onsager},~\bfnm{Lars}\binits{L.}}
(\byear{1944}).
\btitle{Crystal statistics. {I}. {A} two-dimensional model with an
order--disorder transition}.
\bjournal{Phys. Rev. (2)}
\bvolume{65}
\bpages{117--149}.
\bid{mr={0010315}}
\end{barticle}
%
\endbibitem

%b33 ###
\bibitem{parknewman04}
%
\begin{barticle}[vtex]
\bauthor{\bsnm{Park},~\bfnm{Juyong}\binits{J.}} \AND
\bauthor{\bsnm{Newman},~\bfnm{M.~E.~J.}\binits{M.~E.~J.}}
(\byear{2004}).
\btitle{Statistical mechanics of networks}.
\bjournal{Phys. Rev. E (3)}
\bvolume{70}
\bpages{066117--066122}.
\bid{doi={10.1103/PhysRevE.70.066117}, mr={2133807}}
\end{barticle}
%
\endbibitem

%b34 ###
\bibitem{parknewman05}
%
\begin{barticle}[vtex]
\bauthor{\bsnm{Park},~\bfnm{J.}\binits{J.}} \AND
\bauthor{\bsnm{Newman},~\bfnm{M.~E.~J.}\binits{M.~E.~J.}}
(\byear{2005}).
\btitle{Solution for the properties of a clustered network}.
\bjournal{Phys. Rev. E}
\bvolume{72}
\bpages{026136--026137}.
\end{barticle}
%
\endbibitem

%b35 ###
\bibitem{raic04}
%
\begin{barticle}[mr]
\bauthor{\bsnm{Rai{\v{c}}},~\bfnm{Martin}\binits{M.}}
(\byear{2007}).
\btitle{C{LT}-related large deviation bounds based on {S}tein's method}.
\bjournal{Adv. in Appl. Probab.}
\bvolume{39}
\bpages{731--752}.
\bid{doi={10.1239/aap/1189518636}, mr={2357379}}
\end{barticle}
%
\endbibitem

%b23 ###
\bibitem{simon73}
%
\begin{barticle}[mr]
\bauthor{\bsnm{Simon},~\bfnm{Barry}\binits{B.}} \AND
\bauthor{\bsnm{Griffiths},~\bfnm{Robert~B.}\binits{R.~B.}}
(\byear{1973}).
\btitle{The {$(\phi\sp{4})\sb{2}$} field theory as a classical
{I}sing model}.
\bjournal{Comm. Math. Phys.}
\bvolume{33}
\bpages{145--164}.
\bid{mr={0428998}}
\end{barticle}
%
\endbibitem

%b36 ###
\bibitem{stein72}
%
\begin{binproceedings}[mr]
\bauthor{\bsnm{Stein},~\bfnm{Charles}\binits{C.}}
(\byear{1972}).
\btitle{A bound for the error in the normal approximation to the distribution
of a sum of dependent random variables}.
In \bbooktitle{Proc. Sixth Berkeley Symp.
Math. Statist. Probab. Vol. II: Probability Theory}
\bpages{583--602}.
\bpublisher{Univ. California Press}, \baddress{Berkeley, CA}.
\bid{mr={0402873}}
\end{binproceedings}
%
\endbibitem

%b37 ###
\bibitem{stein86}
%
\begin{bbook}[vtex]
\bauthor{\bsnm{Stein},~\bfnm{Charles}\binits{C.}}
(\byear{1986}).
\btitle{Approximate Computation of Expectations}.
%Notes---Monograph Series}
\bpublisher{IMS}, \baddress{Hayward, CA}.
\bid{mr={882007}}
\end{bbook}
%
\endbibitem

%b38 ###
\bibitem{talagrand95}
%
\begin{barticle}[mr]
\bauthor{\bsnm{Talagrand},~\bfnm{Michel}\binits{M.}}
(\byear{1995}).
\btitle{Concentration of measure and isoperimetric inequalities in product
spaces}.
\bjournal{Publ. Math. Inst. Hautes \'Etudes Sci.}
\bvolume{81}
\bpages{73--205}.
\bid{mr={1361756}}
\end{barticle}
%
\endbibitem

%b39 ###
\bibitem{vu01}
%
\begin{barticle}[mr]
\bauthor{\bsnm{Vu},~\bfnm{Van~H.}\binits{V.~H.}}
(\byear{2001}).
\btitle{A large deviation result on the number of small subgraphs of a random
graph}.
\bjournal{Combin. Probab. Comput.}
\bvolume{10}
\bpages{79--94}.
\bid{mr={1827810}}
\end{barticle}
%
\endbibitem

\end{thebibliography}
\end{document}